\documentclass{scrartcl}

\usepackage{natbib}
\usepackage{hyperref}
\usepackage{sidecap}
\usepackage{wrapfig}
\usepackage{lscape,textcomp}
\usepackage[T1]{fontenc}
\usepackage{xparse}
\usepackage{amsmath,amssymb,dsfont,units,booktabs}

\usepackage{amsfonts,import,mathtools}

\usepackage{graphicx, color}
\usepackage{color}
\usepackage[utf8]{inputenc}

\newtheorem{algo}{Algorithm}
\newtheorem{remark}{Remark}

\newtheorem{theorem}{Theorem}

\makeatletter
\newcommand{\opnorm}{\@ifstar\@opnorms\@opnorm}
\newcommand{\@opnorms}[1]{%
  \left|\mkern-1.5mu\left|\mkern-1.5mu\left|
   #1
  \right|\mkern-1.5mu\right|\mkern-1.5mu\right|
}
\newcommand{\@opnorm}[2][]{%
  \mathopen{#1|\mkern-1.5mu#1|\mkern-1.5mu#1|}
  #2
  \mathclose{#1|\mkern-1.5mu#1|\mkern-1.5mu#1|}
}
\makeatother

\DeclareDocumentCommand{\spL}{ s m }{
  \IfBooleanTF #1
  {\left(#2\right)}
  {(#2)}
}
\newcommand\spTS[1]{(\mkern-4mu(#1)\mkern-4mu)}
\newcommand\spLL[1]{(#1)}

\newcommand{\vt}{\mathbf{v}}

\newcommand{\HT}{\mathbf{h}}

\newcommand{\wt}{\mathbf{w}}

\newcommand{\zt}{\mathbf{z}}

\newcommand{\sigmat}{\boldsymbol{\sigma}}
\newcommand{\phit}{\boldsymbol{\phi}}
\newcommand{\psit}{\boldsymbol{\psi}}
\newcommand{\epsilont}{\boldsymbol{\epsilon}}

\newcommand{\taut}{\boldsymbol{\tau}}
\newcommand{\te}{\vec{e}}

\renewcommand{\div}{\operatorname{div}}
\newcommand{\tr}{\operatorname{tr}}
\newcommand{\qt}{\mathbf{q}}

\begin{document}

\title{A goal oriented error estimator and mesh adaptivity for sea ice simulations}  

\author{\underline{Carolin Mehlmann}\thanks{Max-Planck-Institute of Meteorology,
    Bundesstrasse 53, 22176  Hamburg, Germany,
    carolin.mehlmann@mpimet.mpg.de}  
  \and Thomas Richter\thanks{Institute of Analysis and Numerics,
    Otto-von-Guericke University Magdeburg, Universitätsplatz 2,
  39106 Magdeburg, Germany, thomas.richter@ovgu.de}}

\maketitle

\begin{abstract}
  For the first time we introduce an error estimator for the numerical
  approximation of the equations describing the dynamics of sea
  ice. The idea of the estimator is to identify different error contributions coming from
  spatial and temporal discretization as well as from the splitting in time of
  the ice momentum equations from further parts of the coupled
  system. The novelty of the error estimator lies in the consideration of the splitting error, which turns out to be dominant with increasing mesh resolution.  Errors are measured in user specified functional outputs
  like the total sea ice extent. The error estimator is
  based on the dual weighted residual method that asks for the
  solution of an additional dual problem for obtaining sensitivity
  information. Estimated errors can be used to validate the accuracy
  of the solution and, more relevant, to reduce the discretization error by guiding an adaptive
  algorithm that optimally balances the mesh size and the time step
  size to increase the efficiency of the simulation.  
\end{abstract}

\section{Introduction}
We consider the viscous-plastic (VP) sea model, that was introduced by
Hibler in 1979 and which is still one of the most widely used sea ice
rheologies as detailed by~\cite{Stroeve2014}. The model 
includes strong nonlinearities such that solving the sea ice dynamics
at high resolutions is extremely costly and good solvers are under active
research. Mostly, solutions to the VP model are approximated by 
iterating an elastic-viscous-plastic (EVP) modification of the
model that was introduced by~\cite{Hunke1997} and that allows for explicit sub-cycling. Alternatively the VP model is tackled directly with simple Picard iterations as described by \cite{Hibler1979}
or solved with Newton-like 
methods as described by~\cite{Lemieux2010} or ~\cite{MehlmannRichter2016newton}. All approaches are not
satisfactory as they are 
extremely expensive and often are not able to give an accurate solution in
reasonable computational time. It is therefore  of utmost importance to 
reduce the complexity of the computations, e.g. by using coarse meshes
and large time step sizes, as long as this does not deteriorate the accuracy
assumptions.

We derive an error estimator that 
identifies the errors coming from spatial and temporal
discretization. Furthermore, the error estimator allows for a
localization of the error to each element and each time step 
such that local step sizes can be adjusted.
This goal oriented error estimator for the viscous-plastic sea model
is an extension of the dual weighted residual method that was
introduced by~\cite{BeckerRannacher2001}. The aim of the estimator
is to identify discretization errors $J(U)-J(U_{k,h})$ between the unknown exact solution
$U$ and the numerical approximation $U_{k,h}$, where $k$ indicates the
temporal and $h$ the spatial discretization parameter, in functionals
$J(\cdot)$. These functionals can be any measures of interest, e.g. the
average sea 
ice extent in a certain time span
\begin{equation}\label{functional}
  J(U_{k,h}) = \frac{1}{t_2-t_1}\int_{t_1}^{t_2}\int_\Omega
  A_{k,h}(x,y,t)\,\text{d}(x,y)\,\text{d}t. 
\end{equation}
We denote by $A_{k,h}$ the ice concentration (one component of
the solution $U_{k,h}$ which will be introduced later), by $\Omega \subset \mathds{R}^2$ the spatial
domain of interest and by $[t_1,t_2]$ the time span of interest,
e.g. the summer months. The error estimator 
will give approximations to $J(U)-J(U_{k,h})$ which can be attributed
to a spatial error, a temporal error and to a splitting error - coming from
partitioning the system  into momentum equation and balance laws, which is the
standard procedure in sea ice numerics, see~\cite{Lemieux2014}. The
estimation 
of errors in space and time for parabolic problems was discussed by~\cite{Schmich2012}. For the first time, we extend the
application to the VP model and additionally
consider the splitting error. 
\cite{Lipscomb2007} pointed out that decoupling the
system in time can lead to a numerical unstable solution such that a
small time step is required to achieve a stable
approximation. \cite{Lemieux2014} introduced an
implicit-explicit time integration method (IMEX), which resolves this issues
and allows the use of larger time steps. The error estimator will be
able to predict the accuracy implications of this temporal splitting.

Furthermore, it provides information about the spatial and temporal convergence of the approximation of the solution. So far spatial convergence has been analyzed by~\cite{WILLIAMS2018}, where a one dimensional test case has been studied.
The authors observe that the simulated velocity field
depends on the spatial resolution and found that the mean sea ice
drift speed rises by 32\% by increasing resolution from 40 km to 5 km.
The temporal and spatial scaling properties of the mean
deformation rate and the sea ice thickness are studied
by~\cite{Hutter2018}. 

The dual weighted residual estimator by~\cite{BeckerRannacher2001} relies on a variational formulation of the system of partial differential equations and has been introduced for the finite element method. Later on, the estimator has been extended to time dependent problems by using a relation between classical time stepping schemes like the Euler method and temporal Galerkin methods, see~\cite{Schmich2008}, which can be considered as finite elements in time. This similarity is also exploited in this work. The temporal Galerkin approach appears abstract, but it allows for a simple realization in the context of standard time stepping schemes like the backward Euler method and it is necessary in order to formulate the estimator.

Likewise, finite volume methods can be interpreted as discontinuous Galerkin methods that also give direct access to the framework of the dual weighted residual estimator, see~\cite{AfifBergamMghazliVerfuerth2003} or~\cite{ChenGunzburger2014} for an approach that does not rely on this similarity to Galerkin methods and which shows applications in climate modeling. An application of the goal oriented error estimator in finite difference discretizations is less natural, since the variational structure is missing. This however is the basis for the definition of adjoint problems and also of the residual terms that form the estimator. However, low order finite elements are closely related to finite difference methods obtained by numerical quadrature. This similarity is used by~\cite{MeidnerRichter2015} to apply the error estimator to efficient finite difference time stepping schemes and \cite{CollinsEstepTavener2014} exploit the similarity of finite difference schemes with related finite volume and finite element formulations to carry over the idea of the error estimator to spatial finite difference discretizations of conservation laws.

The paper is structured as follows. In Section \ref{sec:comp system}
we start by presenting the sea ice model in strong and variational
formulation which is required for the  Galerkin finite element discretization in space and
time. Further we give details on the partitioned solution approach. In
Section~\ref{sec:error_estimator} we derive the goal oriented error
estimator for the sea ice model and describe its numerical
realization.  We numerically analyse the error estimator in
Section~\ref{sec:num_ana} and conclude in Section~\ref{sec:con}. For
better readability we keep the mathematical formulation as simple as
possible and refer to the literature for  details. Some details on variational formulations are given in the appendix. 


\section{Model Description and Discretization}\label{sec:comp system}

Let $\Omega\subset\mathds{R}^2$ be the spatial domain. We denote  the
time interval of interest by $I=[0,T]$. Sea ice is described by three variables, the sea ice
concentration $A$, the mean sea ice thickness $H$ and the sea ice velocity
$\vt$, such that the complete solution is given by
$U=(\vt,A,H)$.
The VP sea ice model as introduced by~\cite{Hibler1979}
consists of the momentum equation and the balance laws
\begin{equation}\label{eq}
  \begin{aligned}
    \rho_\text{ice}H \big(\partial_t \vt 
    +  f_c \te_r\times (\vt-\vt_\text{ocean})\big)
    &=\div\,\sigmat +
    \taut(\vt),\\
    \partial_t A + \div\,(\vt A) = S_A,\quad 
    \partial_t H &+ \div\,(\vt H) = S_H
  \end{aligned}
\end{equation}
with $0\le H$ and $0\le A\le 1$. 
The forcing term $\taut(\vt)$ models ocean and atmospheric traction
\[
  \taut(\vt) =
  C_\text{ocean}\rho_\text{ocean}\|\vt_\text{ocean}-\vt\|_2(\vt_\text{ocean}-\vt)+
  C_\text{atm}\rho_\text{atm}\|\vt_\text{atm}\|_2\vt_\text{atm},
\]
with the ocean velocity $\vt_\text{ocean}$ and the wind velocity
$\vt_\text{atm}$. 
By $\rho_\text{ice}$ we denote the ice density, by $f_c$ the Coriolis
parameter, by $\te_r$ the radial ($z$-direction) unit vector.
Following~\cite{Coon1980} we have replaced surface height effects
by the approximation $g\nabla \tilde H_g=-f_c\te_r \times
\vt_\text{ocean}\approx 0$. In this paper we focus on the dynamical
part of the sea ice model such that we neglect thermodynamic effects
and set $S_A=0$ and $S_H=0$. 

\begin{table}[t]
  \begin{center}
    \begin{tabular}{l|l|l}
      \toprule
      \text{Parameter} & \text{Definition} & \text{Value}\\
      \midrule
      $\rho_\text{ice}$ & sea ice density &$\unit[900]{kg/m^{3}}$\\
      $ \rho_\text{atm}$& air density & $\unit[1.3]{kg/m^{3}}$\\
      $\rho_\text{ocean}$& water density & $\unit[1026]{kg/m^{3}}$\\
      $C_\text{atm}$&air drag coefficient &$\unit[1.2] \cdot {10^{-3}}$\\
      $C_\text{ocean}$&water drag coefficient& $\unit[5.5] \cdot {10^{-3}}$\\
      $f_c$ &Coriolis parameter&$\unit[1.46] \cdot \unit[10^{-4}]{s^{-1}}$ \\
      $P^{\star}$&ice strength parameter&$\unit[27.5]\cdot \unit[ 10^{3}]{N/m^2}$\\
      $C$&ice concentration parameter&$20$\\
      \bottomrule
    \end{tabular}
 \caption{Physical parameters of the momentum equation\label{Con}.}
  \end{center}
\end{table}

The system of equations~(\ref{eq}) is
closed by Dirichlet conditions $\vt=0$ on the boundary of the domain
and initial conditions $H(0)=H^0$,
$A(0)=A^0$ and $\vt(0)=\vt^0$ for mean ice
thickness, concentration and velocity at time $t=0$. 

Finally, we present the nonlinear viscous-plastic rheology which
relates the stress $\sigmat$ to the strain rate
\[
\dot\epsilont=\frac{1}{2}\Big(\nabla\vt+\nabla\vt^T\Big),\quad
\dot\epsilont':=\dot\epsilont-\frac{1}{2}\operatorname{tr}(\dot\epsilont)I ,
\]
where $\operatorname{tr}(\cdot)$ is the trace. The rheology  is given by
\begin{equation}\label{model:stress}
  \begin{aligned}
    \sigmat &= 2\eta \dot\epsilont' + \zeta \tr(\dot\epsilont)I -
    \frac{P}{2} I,
  \end{aligned} 
\end{equation}
with the viscosities $\eta$ and $\zeta$, given by
$\eta=\frac{1}{4}\zeta$ and
\begin{equation}\label{deltalimit}
  \zeta=\frac{P}{2\Delta(\dot\epsilont)},\quad
  \Delta(\dot\epsilont):= 
  \sqrt{\frac{1}{2} \dot\epsilont':\dot\epsilont'+\tr(\dot\epsilont)^2 +
    \Delta_{min}^2}.
\end{equation}
$\Delta_{min}= 2\cdot 10^{-9}$ is the threshold that describes
the transition between the viscous and the plastic regime. 
The ice strength $P$ in (\ref{model:stress}) is modeled as
\begin{equation}\label{icestrength}
  P(H,A)= P^\star H \exp\big(-C(1-A)\big),
\end{equation}
with the constant $C=20$. All problem parameters are collected in
Table~\ref{Con}.


\subsection{Variational formulation and discretization}\label{sec:var}

The dual weighted residual estimator by~\cite{BeckerRannacher2001}
relies on a variational formulation of the system of partial
differential equations in space and time and on Galerkin
discretizations (like the finite element method) that discretize the
problem by restricting the admissible space for finding the discrete
solution. In our approach we use a linear finite element
discretization in space and a discretization by piecewise constant
functions in time, which also can be considered as a type of finite
element discretization. This time discretization corresponds to the usual
backward Euler method.

To start with, we multiply the three equations~(\ref{eq}) with test
functions $\phit$, $\psi_A$ and $\psi_H$ and integrate in space and
time
\begin{multline}\label{var1}
  \int_0^T\Big(\spLL{\rho_\text{ice}H \partial_t \vt 
    + f_c \te_r\times \vt -\rho_\text{ice}H f_c\te_r
    \times \vt_\text{ocean},\phi}\\
  \quad  +\spLL{\taut(\vt),\phit}
  +\spLL{\sigmat,\nabla\phit} 
  + \spLL{\partial_t A + \div\,(\vt A),\psi_A}\\
  + \spLL{\partial_t H + \div\,(\vt H),\psi_H}\Big)\,\text{d}t=0. 
\end{multline}
By $\spLL{f,g}\coloneqq \int_\Omega
f(x)g(x)\,\text{d}x$
we denote the usual  $L^2$-inner product. Apart from the
additional integration in time, this is the usual variational
formulation for finite element discretizations as used
in~\cite{Danilov2015}. 

In space, the discretization is briefly described: We define a  
conforming finite element space $V_h$ for velocity $\vt\in V_h^2$, ice
concentration and mean sea ice thickness $A,H\in V_h$. In our 
implementation we use the space of piecewise bi-linear functions
defined on a quadrilateral mesh $\Omega_h$ of the domain
$\Omega$. \cite{Danilov2015},
  \cite{Dansereau2016}, \cite{Rampal2016} consider finite element
  discretization  on triangular meshes. The error estimator presented
  in the following section directly transfers to such discretizations,
  as it is shown in~\cite{CarpioPrietoBermejo2013}.

To discretize (\ref{var1}) in time we split the interval $I=[0,T]$
into equidistant\footnote{For simplicity we assume that
  $k=t_n-t_{n-1}$ is constant for all the time steps. The extension to
  varying time step sizes is discussed in
  literature, see ~\cite{Schmich2008}.} discrete steps
$0=t_0<t_1<\cdots<t_N=T$ with $k\coloneqq t_n-t_{n-1}$. By
  discretizing with piecewise constant functions in time and by
  introducing appropriate scalar products that penalize the jumps of
  $\vt$ and $H,A$ at each time step $t_n$ (compare Figure \ref{Def:spruenge} in the Appendix for a visualization of the jumps), we obtain the
  standard implicit Euler discretization of the finite element
  formulation. For $n=1,2,\dots,N$ we solve
\begin{equation}\label{eulerstep1}
  \begin{aligned}
    \spLL{\rho_\text{ice}H_{n}(\vt_{n}-\vt_{n-1}),\phit_h}
    +k\spLL{
      \rho_\text{ice}H_nf_c \te_r\times (\vt_{n} -\vt_\text{ocean}),\phit_h}\qquad &\\ 
    +k\spLL{\taut(\vt_{n}),\phit_h}
    +\spLL{\sigmat(\vt_{n},A_{n},H_{n}),\nabla\phit_h}&=0 \\
    \spLL{A_{n}-A_{n-1},\psi_A}+
    k \spLL{\div\,(\vt_{n} A_{n}),\psi_{A}}&=0\\
    \qquad\spLL{H_{n}-H_{n-1},\psi_H}+
    k\spLL{ \div\,(\vt_{n} H_{n}),\psi_{H}}&=0\\
  \end{aligned}
\end{equation}
Division by the step size $k$ reveals the classical backward Euler
scheme which is standard in sea ice dynamics as described
by~\cite{Lemieux2014}. The discrete functions
$\vt_n,A_n,H_n,\phit_h,\psi_A,\psi_H$ are all piecewise linear finite
elements in space.
\begin{remark}
  The transport equations for $A$ and $H$ are under the constraints
  $0\le H$ and $0\le A\le 1$ which is not easily accessible for a variational formulation. We will therefore realize this constraint weakly by introduction of the following right hand side
  \[
  \spLL{A_{n}-A_{n-1},\psi_A}+
  k \spLL{\div\,(\vt_{n} A_{n}),\psi_{A}}
  =\spLL{\min\{0,1-A_n\},\psi_A} 
  \]
  that only gets active if $A_n>1$ and that will then force $A_n$
  below one.
\end{remark}

To present the dual weighted residual error estimator we introduce
the notation of the Galerkin formulation which is equivalent to the
implicit Euler formulation given in~(\ref{eulerstep1})
\begin{multline}\label{var_disc}
  \sum_{n=1}^N 
  \int_{t_{n-1}}^{t_n}\Big(\spLL{ \rho_\text{ice} H\big(\partial_t
    \vt+f_c\te_r\times (\vt-\vt_\text{ocean})\big),\phit_h}\\
  \quad\qquad+\spLL{\sigmat(\vt,A,H),\nabla\phit_h}
  +\spLL{\taut(\vt),\phit_h}\\
  \quad\qquad+\spLL{\partial_t A+\div\,(\vt A),\psi_A}
  +\spLL{\partial_t H+\div\,(\vt H),\psi_H}\Big)  \,\text{d}t\\
  \qquad
  +\spLL{\rho_\text{ice}H(t_{n-1})^+[\vt]_{n-1},\phit(t_{n-1})^+}\\
  \qquad
  +\spLL{[A]_{n-1},\psi_A(t_{n-1})^+}
  +\spLL{[H]_{n-1},\psi_H(t_{n-1})^+}=0, 
\end{multline}
where the discrete function space, which contains piecewise constant functions in time and linear finite elements in space, is denoted by  $V_{k,h}$. As the discrete functions of  $V_{k,h}$ are discontinuous at each step $t_n$, we 
denote by $H(t_n)^+$ and $H(t_n)^-$ their values from the right and
the left and by $[H]_n=H(t_n)^+-H(t_n)^-$ the jump of
discontinuity. We also refer to Appendix~\ref{app:gd} for details
and to~\cite{Thome1997} for a comprehensive background on temporal Galerkin discretizations.
The real solution is continuous in time and it holds $[\vt]_n=0$ and 
$[A]_n=[H]_n=0$ such that the true solution to~(\ref{eq})
or~(\ref{var1}) is also a solution to this discrete
formulation~(\ref{var_disc}).
The beauty of the Galerkin approach lies in the
presence of one single problem formulation~(\ref{var_disc}) that is 
equivalent to the original problem~(\ref{eq}) if arbitrary functions
are allowed for solution $\vt,\HT$ and for testfunctions $\phit,\psit$ and
that is equivalent to the Euler  / finite element
discretization~(\ref{eulerstep1}) if solution and test functions are
restricted to piecewise constants in time and piecewise linear
functions in space. We solve the problem with the established and
efficient Euler scheme but we fall back to the Galerkin scheme in an
integral formulation, when it comes to estimating the error, see
Section~\ref{Goal:IMEX}. 

%

To shorten notation we
combine $U=(\vt,\HT)$ with $\HT=(A,H)$ and $\Phi=(\phit,\psit)$ with
$\psit=(\psi_A,\psi_H)$ and we assume that these functions come from
function spaces $U\in {\cal X}:={\cal V}\times {\cal V}^A\times {\cal
  V}^H$ and $\Phi\in {\cal X}$. The exact notation of all function
spaces is introduced in~\cite{MehlmannRichter2016mg}. Then, the
variational formulation 
(\ref{var_disc})  can
written in an abstract notation by introducing the form $B(U)(\Phi)$
which simply collects all the integrals and jumps from~(\ref{var_disc})
\begin{equation}\label{varform1}
  U\in {\cal X}\quad B(U)(\Phi) = 0\quad\forall \Phi\in {\cal X}.
\end{equation}
The discrete solution
$U_{k,h}=(\vt_{k,h},A_{k,h},H_{k,h})\in X_{k,h}:=V_{k,h}\times V_{k,h}^A\times
V_{k,h}^H$ is given by restricting~(\ref{varform1}) to the finite 
dimensional discrete space $X_{k,h}$
\begin{equation}\label{vardisc1}
  U_{k,h}\in X_{k,h}\quad B(U_{k,h})(\Phi_{k,h}) = 0\quad\forall
  \Phi_{k,h}\in X_{k,h}. 
\end{equation}

\subsection{Partitioned solution approach}\label{sec:part}

The discrete formulation~(\ref{vardisc1}) naturally splits into
time-steps $t_{n-1}\to t_n$ as shown in~(\ref{eulerstep1}). The three
components velocity $\vt$, ice  concentration $A$ and mean ice thickness
$H$ however are coupled. It is standard 
to apply a partitioned solution approach in every time step, either by
first solving the momentum equation for the sea ice velocity followed
by the balance laws, or vice versa, see~\cite{Lemieux2014}. We
first solve the momentum equation and replace all appearances
of the ice concentration $A_n$ and mean ice thickness $H_n$ in the
momentum equation of~(\ref{eulerstep1}) by the previous
approximations $H_{n-1}$ and $A_{n-1}$.

To realize this decoupling within the Galerkin formulation~(\ref{var_disc}), we introduce 
the projection operator ${\cal S}:{\cal X}\to {\cal X}$ that
projects $A\in {\cal V}^A$ (or $H$, respectively) on the interval
$I_n=(t_{n-1},t_n]$ onto the value $H(t_{n-1})^-$. For discrete
functions this corresponds to $S(H_n) = H_{n-1}$ and
$S(A_n)=A_{n-1}$. 
This calls for a slight modification of
the variational formulation $B(U)(\Phi)=0$ in (\ref{varform1}), namely the introduction of
the projection operator ${\cal S}$ in the momentum part
\begin{multline}\label{varsplit}
  B_s(U)(\Phi)=
  \sum_{n=1}^N 
  \int\limits_{t_{n-1}}^{t_n}\Big(\spLL{ \rho_\text{ice} {\cal S}(H)\big(\partial_t
    \vt+f_c\te_r\times (\vt-\vt_\text{ocean})\big),\phit_h}\\
  \quad\qquad+\spLL{\sigmat(\vt,{\cal S}(A),{\cal S}(H)),\nabla\phit_h}
  +\spLL{\taut(\vt),\phit_h}\Big)\text{d}t\\
  \qquad
  +\spLL{\rho_\text{ice}{\cal S}(H)(t_{n-1})^+[\vt]_{n-1},\phit(t_{n-1})^+}+\dots
\end{multline}
where the dots denote the equations for $A$ and $H$, which are not
changed in comparison to~(\ref{var_disc}). 
Once again, we indicate this variational form only for the formulation and evaluation
of the error estimator, the solution itself is computed by the
backward Euler scheme~(\ref{eulerstep1}) by replacing $A_n,H_n$ in the
momentum equations by the previous approximations $A_{n-1},H_{n-1}$.
Since the discretization is based on~(\ref{varsplit}), but the
exact solution is given by~(\ref{var_disc}), the resulting
discretization is called a non-consistent Galerkin formulation, see Appendix \ref{app:go} for details.

\section{Goal oriented error estimation}\label{sec:error_estimator}

 In this section, we derive a \emph{goal oriented error estimator} for
 \emph{partitioned solution approaches}. 
The new error estimator will be based on concepts of the \emph{dual
  weighted residual} (DWR)  method introduced
by~\cite{BeckerRannacher2001}. The DWR estimator can be easily applied
to 
all problems given in a variational Galerkin formulation and it has
been applied to various problems for error estimation in space and
time such as fluid dynamics~(\cite{Schmich2012}) or fluid-structure
interactions \cite[Chapter 8]{Richter2017}. 

The novel aspect of our approach is to properly include the
error that stems from using a partitioned solution approach. 
Splitting will result in a non-consistent variational formulation,
the analytical problem and the discrete problem do not match anymore. Measuring  this splitting error will call for additional 
effort. This will be discussed in
the following section and also in Appendix~\ref{app:go}. Further, in Section~\ref{sec:num1}, we discuss the various
parts that make up the error estimator. It will turn out that
including this splitting error is essential as its share in the
overall error can be dominant.

\subsection{The goal oriented error estimator for
  partitioned solution approaches}  
\label{Goal:IMEX}
The dual weighted residual estimator measures the
discretization error with respect to a goal functional $J(\cdot)$
like~(\ref{functional}), measuring the average sea ice extent. It
  is formulated as an optimization 
problem: we minimize  $J(U)$ under the constraint that $U$ solves
the sea ice problem. To tackle this optimization problem we introduce
the Lagrangian  
\begin{equation}\label{LAG1}
L(U,Z) \coloneqq J(U) - B(U)(Z),
\end{equation}
where $Z$ takes the role of the Lagrange multiplier. Due to the splitting approach in time  the discrete solution solves (\ref{varsplit}), but not (\ref{var_disc}). Thus,  we introduce a second Lagrangian. 
\begin{equation}\label{LAGs1}
  L_s(U_{k,h},Z_{k,h})\coloneqq J(U_{k,h})-B_s(U_{k,h})(Z_{k,h}),
\end{equation}
which is based on the projection operator ${\cal S}$. This
differentiates the error estimate from the standard dual weighted
residual method, where it is sufficient to use only one Lagrangian.

Since $B(U)(Z)=0$ for the true solution and $B_s(U_{k,h})(Z_{k,h})=0$
for the discrete solution, we obtain the nonlinear error identity 
\begin{equation}\label{DWR0}
  J(U)-J(U_{k,h}) = L(U,Z) - L_s(U_{k,h},Z_{k,h}),
\end{equation}

We only sketch the derivation and refer
to~\cite{Mehlmann2019} for details. The non-consistency of the 
variational formulation coming from the splitting approach is
incorporated by introducing $\pm L(U_{k,h},Z_{k,h})=0$ and by
separating the error influences into the Galerkin error (from
discretization in space and time) and the  splitting error (from partitioned time
stepping)
\[
J(U)-J(U_{k,h}) =
\underbrace{L(U,Z)-L(U_{k,h},Z_{k,h})}_\text{Galerkin}
+
\underbrace{L(U_{k,h})(Z_{k,h})-L_s(U_{k,h},Z_{k,h})}_\text{splitting}. 
\]
The estimation of the Galerkin part is the standard procedure of the dual weighted Galerkin method.  We reformulate 
\begin{equation}\label{DWR11}
  L(U,Z)-L(U_{k,h}, Z_{k,h})=
  \int_0^1 \frac{d}{ds}L\big(U_{k,h}+s(U-U_{k,h}),
  Z_{k,h}+s(Z-Z_{k,h})\big)\,\text{d}s
\end{equation}
  and define the directional derivative of $L(U,Z)$ in an arbitrary direction $(\psit,\phit)$ as
  \begin{equation}
    L'(U,Z)(\psit,\phit)\coloneqq \frac{d}{ds} L(U+s\psit,Z+s\phit)\Big|_{s=0}. 
  \end{equation}
  To shorten the notation we combine $X\coloneqq (U,Z)$ and $X_{k,h}\coloneqq (U_{k,h},Z_{k,h})$ and approximate the integral in (\ref{DWR11}) with the trapezoidal rule
  \begin{multline}\label{DWR21}
    L(U,Z)-L(U_{k,h}, Z_{k,h}) = L(X)-L(X_{k,h})\\
    =  \int_0^1 L'\big(X_{k,h}+s(X-X_{k,h})\big)\big(X-X_{k,h}\big)\,\text{d}s\\ 
    = \frac{1}{2}L'(X)(X-X_{k,h}) + \frac{1}{2}L'(X_{k,h})(X-X_{k,h})\\
    + \frac{1}{2}\int _0^1s(1-s)L'''\big(X_{k,h}+s(X-X_{k,h})\big)
      \big(X-X_{k,h}\big)\,\text{d}s,
  \end{multline}
  where we denote by $L'''(X)(\Psi)$ the third directional derivative of $L(X)=L(U,Z)$ in direction $\Psi$. See~\cite[Sec. 9.2.2]{QuarteroniSaccoSaaleri2007}  for a derivation of the trapecoidal rule's error formula. The remainder is of third order in the error $X-X_{k,h}=(U-U_{k,h},Z-Z_{k,h})$ and omitted in practical application. 
If we consider the definition of the
Lagrangian~(\ref{LAG1}) the derivatives $L'$ take the form
\[
  L'(U,Z)(\delta U,\delta Z) =
  J'(U)(\delta U) - B'(U)(\delta U,Z)
  -B(U,\delta Z).
\]
Details are given in Section~\ref{sec:realization}. 
To proceed with the nonlinear error identity~(\ref{DWR21}) we now
define $Z\in {\cal X}$  as the solution to
the \emph{linearized adjoint problem}
\begin{equation}\label{DWR3}
  \begin{aligned}
    B'(U)(\Psi,Z) &= J'(U)(\Psi)& \Psi&\in {\cal
      X}
  \end{aligned}
\end{equation}
Analogously one differentiates $L_s'$  and gets
\begin{equation}\label{DWR3a}
  \begin{aligned}
    B_s'(U)(\Psi,Z) &= J'(U)(\Psi)& \Psi&\in {\cal
      X}\\
    B_s'(U_{k,h})(\Psi,Z_{k,h}) &=
    J'(U_{k,h})(\Psi)& \Psi&\in  X_{k,h},
  \end{aligned}
\end{equation}
where we define  $Z_{k,h}\in {\cal X}$  as the solution to discretized
\emph{linearized adjoint problem}.
\begin{remark}[Adjoint solution]
  The solution $Z$ to the linearized adjoint problem~(\ref{DWR3})
  indicates the sensitivity of the error functional $J(U)$ with
  respect to the variations in the solution $U$. The adjoint solution runs backward in time and the direction of transport is reversed. The use of adjoint
  solutions is standard in the area of constraint optimization
  problems, where the adjoint solution takes the role of the
  Lagrange multiplier. In sea ice models (in general in all parts of
  climate models), adjoint equations play an eminent role in
  variational data assimilation, which can be
  considered as gradient based calibration of the model with respect
  to measurement data. The adjoint problems required for the process
  of error estimation result from the same linearized equations,
  with the right hand side given by the goal functional
  $J(\cdot)$. Some climate models like \emph{MITgcm} (see \cite{adjointmitGCM1999,HeimbachHillGiering2005}) or \emph{MRI.COM} (see \cite{mricom}) support adjoint equations for coupled ice, ocean and atmosphere simulation. This will simplify the realization of the error estimator in climate models. 
\end{remark}
Equation (\ref{DWR21}) still depends on the discretization errors $U-U_{k,h}$ and $Z-Z_{k,h}$ which are unknown. However, these errors can be replaced by interpolation errors $U-I_{k,h}U$ and $Z-I_{k,h}Z$ which can be approximated by local reconstructions. For details we refer to Appendix \ref{app:go}, see also Remark~\ref{rem:weights}.

The splitting part is derived by using the definitions of the
Lagrangians in~(\ref{LAG1}) and (\ref{LAGs1})
\begin{equation}\label{Lsplit}
  L(U_{k,h},Z_{k,h})-L_s(U_{k,h},Z_{k,h})
  =\beta(U_{k,h})(Z_{k,h}),
\end{equation}
where the \emph{splitting error} is given by the difference between
original variational form~(\ref{var_disc}) and splitting form~(\ref{varsplit})
\begin{equation}\label{consistency}
  \beta(U)(\Phi):=B_s(U)(\Phi)-B(U)(\Phi).
\end{equation}

This error contribution can be evaluated since it only depends on
quantities that are available, namely the primal and dual discrete
solutions. We summarize: 

\begin{theorem}[DWR estimator for partitioned solution
    schemes] 
  \label{Theo:2}
  Let $U,Z\in {\cal X}$ and $U_{k,h}, X_{k,h}\in  {\cal X}_{k,h}$ be
  primal and dual solutions to
  \[
  \begin{aligned}
    B(U)(\Phi) &= 0 & \forall \Phi&\in {\cal X}\\
    B'(U)(\Psi,Z)&= J'(U)(\Psi)& \forall \Psi&\in {\cal X}\\
    B_s(U_{k,h})(\Phi_{k,h}) &= 0& \forall \Phi_{k,h}&\in X_{k,h}\\
    B_s'(U_{k,h})(\Psi_{k,h},Z_{k,h})&= J'(U_{k,h})(\Psi_{k,h})
    & \forall \Psi_{k,h}&\in X_{k,h}.
  \end{aligned}
  \]
  Then, it holds that 
  {\footnotesize
    \begin{multline}\label{error_est}
      J(U)-J(U_{k,h})= {\mathcal{R}}(x_{k,h},e)
      -\frac{1}{2} B_s(U_{k,h})(Z-I_{k,h} Z)\\
      +\frac{1}{2}\Big\{
      J'(U_{k,h})(U-I_{k,h} U)-
      B'_s(U_{k,h})(U-I_{k,h}U,Z_{k,h})
      \Big\}\\
      + \frac{1}{2}\Big\{\beta(U_k)(Z+Z_{k})
      + \beta'(U_k)(U_k,Z-Z_k)\Big\}
  \end{multline}}
  where ${\mathcal{R}}(x_{k,h},e)$ is given in (\ref{DWR21}) and with the primal and dual splitting errors
  \[
  \begin{aligned}
    \beta(U)(\Phi)&:=B_s(U)(\Phi)-B(U)(\Phi),\\
    \beta'(U)(\Psi,Z)&:=B_s'(U)(\Psi,Z)-B'(U)(\Psi,Z). 
  \end{aligned}
  \]
 By $I_{k,h}:{\cal X}\to X_{k,h}$ we denote an interpolation to the
  space-time domain, by $U_k$ and $Z_k$ we denote semidiscrete
  solutions which are discretized in time only. 
\end{theorem}

\begin{remark}[Weights]\label{rem:weights}
  The error estimator~(\ref{error_est}) depends on the unknown solutions $U$ and $Z$ but
  also on semidiscrete solutions $U_k$ and $Z_k$, which are still
  continuous in space.  All these terms must be approximated by
  suitable reconstruction techniques in order to evaluate the error
  estimator. In general, the reconstruction is realized by a postprocessing mechanism: The discrete solution is reinterpreted as a solution of double polynomial degree (linear in time instead of constant, quadratic in space instead of linear). In time, the discontinuity is resolved and in space we combine adjacent element to form the quadratic function. See Figure~\ref{fig:recon} in the appendix for an illustration of this reconstruction process.
  In time, we combine two intervals
  and reconstruct a linear function by
  connecting $U_{n-1}$ in $t_{n-1}$ with $U_n$ at $t_n$, i.e.
  \begin{equation}\label{recontime}
    i_{k}^{(1)}U_{k,h}\Big|_{(t_{n-1},t_n]}
      = U_{n-1} + \frac{t-t_{n-1}}{k}(U_n-U_{n-1}). 
  \end{equation}
  Then, we approximate $U_k\approx i_{k}^{(1)} U_{k,h}$. In space a
  similar procedure is done by combining the piecewise linear function
  $U_{k,h}$ on four adjacent quadrilaterals to one quadratic
  function. We refer to~\cite{RichterWick2015} and~\cite{Mehlmann2019}
  for details. The notation $i_{k}^{(1)}$ means: interpolation to
  linear  ($i^{(1)}$) functions on the same mesh $i_{k}$. Correspondingly, $i_{2h}^{(2)}$ stands for the interpolation
  to the space of quadratic functions $i^{(2)}$ on the space with double
  spatial mesh spacing $i_{(2h)}$. 
\end{remark}

\subsubsection{Decomposing the error estimator}\label{sec:decomp}

One application of the error estimator is to identify different
contributions to the overall error, namely the error coming from
discretization in space $\eta_{h}$, from the discretization in time $\eta_{k}$  and from the
splitting $\eta_\beta$. This information can help to optimally balance the
discretization, e.g. by avoiding excessive refinement (in space or
time) or by avoiding (or by applying)
a more costly implicit-explicit integration scheme (see~\cite{Lemieux2014}) to avoid the 
splitting error.

The structure of the error estimator~(\ref{error_est}) consists of
residuals weighted by primal, $U-I_{k,h}U$, and dual, $Z-I_{k,h}Z$, interpolation errors and by two 
terms, $\beta$ and $\beta'$, which measure the splitting error $\eta_\beta$. The residual primal and dual weighted residuals refer to
the discretization error in space $\eta_{h}$ and in time $\eta_{k}$. An allocation of this
combined space-time error can be achieved by introducing intermediate
interpolations, which we discuss for the primal residual term $B_s(U_{k,h})(Z-I_{k,h}Z)$, the first term in~(\ref{error_est}). We introduce $\pm I_k Z$, an
interpolation into the space of functions that are piecewise constant
in time but still non-discrete in space
\begin{equation}\label{err:divide}
  B_s(U_{k,h})(Z-I_{k,h}Z) =
  \underbrace{B_s(U_{k,h})(Z-I_{k}Z)}_\text{time}+
  \underbrace{B_s(U_{k,h})(I_k Z-I_{k,h}Z)}_\text{space}.
\end{equation}
We can split the residual into two separate parts, as the dependency
on the weight (which takes the role of the test function) is always
linear.

Naturally, the interpolation $I_k Z$ is not available. However, we can
approximate the interpolation errors by the reconstruction operator
that have been mentioned in Remark~\ref{rem:weights}. To be precise,
the two terms in~(\ref{err:divide}) are approximated by
\[
  B_s(U_{k,h})(Z-I_{k,h}Z)\approx 
  B_s(U_{k,h})(i_{k}^{(1)} Z_{k,h}-Z_{k,h})+
  B_s(U_{k,h})(i_{2h}^{(2)} Z_{k,h}-Z_{k,h}).
\]
The philosophy is simple: for estimating the error in time, we compare
the discrete solution $Z_{k,h}$ with its higher order reconstruction
in time $i_{k}^{(1)}Z_{k,h}$, the spatial error is estimated by
considering the spatial reconstruction operator only. All further
residual terms in~(\ref{error_est}) are handled in the same way.

\subsection{Realization for sea ice dynamics}\label{sec:realization}  
The standard feedback-approach for running adaptive simulations based
on the DWR method is as follows:
\begin{algo}\label{algo:feedback}
  Let $0=t_0<t_1<\cdots<t_N=T$ be the initial time mesh, $\Omega_h$
  the initial spatial mesh. Let $V_{k,h}$ be the resulting
  space-time function space. 
  \begin{enumerate}
  \item Solve the primal problem $U_{k,h}\in V_{k,h}$ 
  \item Solve the dual problem $Z_{k,h}\in V_{k,h}$
  \item Approximate the residual weights $U-I_{k,h}U$, 
    $Z-I_{k,h}z$ 
  \item Evaluate the error estimator~(\ref{error_est})
  \item Stop, if $|J(U)-J(U_{k,h})|$ is sufficiently small
  \item Otherwise use the error estimator to adaptively refine the
    spatial and temporal discretization and restart with a finer
    function space ${\cal V}_{k',h'}$ on the refined mesh
    $\Omega_{h'}$. 
  \end{enumerate}
\end{algo}
In Step 3, the approximation of the residual weights is the delicate part of
the error estimator. Here, we must replace the unknown exact solutions
$U$ and $Z$. For heuristic  approaches we refer
to~\cite{BeckerRannacher2001}, or in particular~\cite{Schmich2008}
or~\cite{MeidnerRichter2014} for time discretizations.

Application of the DWR method will always require some numerical
overhead, mainly by the computation of the auxiliary dual problem. It
turns out that the dual problem is always a linear problem, also in
the case of the fully nonlinear sea ice problem. In the following
sections we describe the steps that are required for applying the DWR
estimator to the sea ice model in a partitioned solution framework. We
refer to~\cite{Mehlmann2019} for the full derivation.  

In order to apply the adaptive feedback loop presented in 
Algorithm~\ref{algo:feedback}, we will first
indicate the exact discrete formulations for solving the primal and
dual problems. These are not written as a discontinuous Galerkin  formulation
but in the form of the Euler method. Afterwards we give some further
remarks on the evaluation of the error estimator. 

\begin{algo}[Primal  sea ice problem] \label{algo:prim}
  Let $\vt_0$ and $\HT_0=(A_0,H_0)$ be the initial solutions at time
  $t=0$. Iterate for $n=1, ...,N$ 
  \begin{enumerate}
  \item Solve for the velocity $\vt_n\in V_h$
    \begin{multline*}
      k^{-1}\spLL{\rho_\text{ice}H_{n-1}(\vt_n-\vt_{n-1}),\phi}_\Omega 
      \\ + 
      \spLL{\rho_\text{ice}H_{n-1} f_c \te_r\times (\vt_n-\vt_\text{ocean}) + \taut(\vt_n),\phit}_\Omega\\ 
      +\spLL{\sigmat(\vt_n,H_{n-1},A_{n-1}),\nabla\phit}_\Omega=0
      \quad\forall \phit\in V_h
    \end{multline*}
  \item Solve the transport equations for $A_n\in V^A_h$ and $H_n\in V^H_h$
    \[
    \begin{aligned}
      k^{-1}\spLL{A_n-A_{n-1},\psi_A}_\Omega
      + \spLL{\div(\vt_nA_n),\psi_A}_{\Omega}
      &=(\min\{0,1-A_n\}, \psi_A)&
      \quad\forall \psi_A&\in V^A_h
      ,\\
      k^{-1}\spLL{H_n-H_{n-1},\psi_H}_\Omega
      +\spLL{\div(\vt_nH_n),\psi_H}_{\Omega}
      &=0\quad&
      \forall \psi_H&\in V^H_h
    \end{aligned}
    \]
  \end{enumerate}
\end{algo}

To derive the dual sea ice model defined in Theorem~\ref{Theo:2}, we must differentiate the form $B_s(U)(\Phi)$, which is described in~(\ref{varsplit}), in the direction of the solution $U=(\vt,\HT)$. The dual solution $Z=(\wt,\qt)$ replaces the test function and the new test function $\Psi=(\psit,\psi_H,\psi_A)$ is the argument of the directional derivative. The complete derivative of the variational formulation is presented in~\cite[p. 131]{Mehlmann2019}. For (rather) simple equations like the viscous plastic sea ice model it is most efficient to implement the dual solution (i.e. the adjoint of the linearization) analytically. The additional effort is minimal if a Newton scheme is employed for the solution of the forward problem since the adjoint system matrix is just the transposed of the Jacobian. In more complex coupled models and in particular, if variants of models are considered, the tangent problem (the Jacobian) and the adjoint can be efficiently generated by automatic differentiation, see~\cite{GriewankWalther2008}. This is approach is implemented in \emph{MITgcm}, see~\cite{HeimbachHillGiering2005}.

Most characteristic for the dual problem is the reversal of the time
direction, the problem runs \emph{backward in time}.
%
This reversal of direction also carries over to the splitting
scheme. While the primal iteration first solves the momentum equation,
the dual problem naturally results in first solving the (dual)
transport problems, followed by the (dual) momentum equation.

\begin{algo}[Partitioned solution approach for the dual
    system]\label{algo:dual} 
  Let $\vt_n$ and $\HT_n$  for $n=0,\dots,N$, be the discrete solution
  of the primal problem. We set $\vt_{N+1}:=0$ and $\zt_{N+1}:=0$ and
  iterate backward in time from $n=N$ to $n=1$
  \begin{enumerate}
  \item Solve the dual transport equations for
    $\qt_n=(\qt_{A,n},\qt_{H,n})\in V_h^2$ 
    \begin{multline*}
      k^{-1}\spLL{\qt_{A,n}-\qt_{A,n+1},\psi_A}_\Omega
      -\spLL{\vt_n\cdot\nabla\qt_{A,n},\psi_A}_\Omega\\
      \qquad+\sigmat'_{A}(\vt_{n+1},H_n,A_n)(\psi_A,\zt_{n+1})
      = J'_{A}(U_n)(\psi_A)
      -\spLL{\chi_{A_n>1},\psi_A}_\Omega
    \end{multline*}
    and
    \begin{multline*}
      k^{-1}\spLL{\qt_{H,n}-\qt_{H,n+1},\psi_H}_\Omega
      -\spLL{\vt_n\cdot\nabla\qt_{H,n},\psi_H}_\Omega\\
      +\spLL{f_c\rho_\text{ice}\zt_{n+1} \te_r\times
        (\vt_{n+1}-\vt_\text{ocean}(t_{n+1})), \psi_H }_{\Omega}\\
      +\sigmat'_{A}(\vt_{n+1},H_n,A_n)(\psi_A,\zt_{n+1})\\
      +k^{-1}\rho_\text{ice}\spLL{(\vt_{n+1}-\vt_n)\zt_{n+1},\psi_H}
      = J'_{H}(U_n)(\psi_H)
    \end{multline*}
    for all $\psi_A,\psi_H\in V_h$.     
  \item Solve the dual momentum equation for $\zt_n\in V_h^2$
    \begin{multline*}
      k^{-1}\spLL{\rho_\text{ice}(H_{n-1}\zt_n-H_n\zt_{n+1}),\phit}_\Omega
      + \spLL{f_c\rho_\text{ice}H_{n-1}\te_r\times \phit,\zt_n}_\Omega\\
      +\spLL{\taut'(\vt_n)(\phit),\zt_n}_\Omega
      +
      \spLL{\sigmat'_{\vt}(\vt_n,H_{n-1},A_{n-1})(\phit,\zt_n)}_\Omega\\
      -\spLL{H_n\nabla \qt_{H,n},\phit}_\Omega
      -\spLL{A_n\nabla \qt_{A,n},\phit}_\Omega=
      J'_\vt(U_n)(\phit)
    \end{multline*}
  \end{enumerate}
  for all $\phit\in V_h^2$. 
  By $\chi_{A>1}(x)$ we denote the characteristic function satisfying
  $\chi_{A>1}(x)=1$ for $A(x)\ge 1$ and $\chi_{A>1}(x)=0$ for
  $A(x)<1$. By $\sigmat'$ we denote the derivatives of the stress
  tensor with respect to $A$, $H$ or $\vt$, by $\taut'$ the
  derivative of the forcing and by $J'$ that of the
  functional. These terms are detailed in~\cite{Mehlmann2019}. 
\end{algo}

\begin{remark}[Dual problem]
  The complexity of the dual system appears immense. However, the
  dual equation is linear such that the solution of each time step
  is very simple and comparable to one single Picard iterations of the forward problem, see~\cite{Lemieux2009}.

  If the sea ice system is linearized by a Newton method, it turns
  out that the dual system matrix is just the transposed of the
  Jacobian. It is hence not necessary to implement the rather
  complicated form of the equations in
  Algorithm~\ref{algo:dual}. Instead, it is sufficient to assemble the
  Newton Jacobian and take its transpose. 
\end{remark}

Primal and dual problem in Algorithm~\ref{algo:prim} and~\ref{algo:dual} are given in the classical Euler time stepping formulation. To evaluate the error estimator~(\ref{error_est}) we must employ the equivalent variational formulations of the discrete forms~(\ref{vardisc1}) and~(\ref{DWR3a}), since the error estimator requires the testing with higher order reconstructions of the weights $i_{k}^{(1)}i_{2h}^{(2)}Z_{k,h}-Z_{k,h}$ and  $i_{k}^{(1)}i_{2h}^{(2)}U_{k,h}-U_{k,h}$, compare Remark~\ref{rem:weights}. Using the variational form is necessary since the equivalence to the Euler scheme only holds for piecewise constant trial and test functions but not for the reconstructed, piecewise linear weights. To approximate all integrals in~(\ref{var_disc}) and the corresponding adjoint form with sufficient accuracy we employ the midpoint rule on each time step. To give an example: the forcing term $(\taut(\vt_n),\phit)_\Omega$ in Algorithm~\ref{algo:prim} corresponds to the space time integral
  $\int_{t_{n-1}}^{t_n} (\taut(\vt),\phit)\,\text{d}t$
  in~(\ref{var_disc}) and with the reconstructed weight $\phit=i_{k}^{(1)}\zt_{k,h}-\zt_{k,h}$ the term within the error estimator is approximated by
\[
\spTS{\taut(\vt_{k,h}),i_{k}^{(1)}\zt_{k,h}-\zt_{k,h}}
=
\frac{k}{2}  \sum_{n=1}^N\big(
\taut(\vt_n),\big(\zt_{n-1}-\zt_n\big)\big)_\Omega,
\]
where we used that $i_{k}^{(1)}\zt_{k,h}-\zt_{k,h}=\frac{1}{2}(\zt_{n-1}-\zt_n)$ in the midpoint of the interval $(t_{n-1},t_n)$, compare Remark~\ref{rem:weights} and~(\ref{recontime}). 
%


\section{Numerical examples}\label{sec:num_ana}

Usually, a posteriori error estimators are used for two objectives: to
compute an approximation with a certain accuracy as stopping criteria for the
simulation, and, to adaptively control the discretization parameters,
namely the mesh size and the time step size.

The first goal is not realistic in sea ice simulations. Uncertainties
from measurement and from model inaccuracies are so large that
quantitative error measures are not available. Furthermore,
large scale simulations are computationally extremely
challenging. Mostly there is little room for using finer and finer
meshes. 
However, the described analysis of the different error contributions
is of great computational importance as it allows to optimally balance
all error contributions to avoid excessive over-refinement in space or
in time. The estimator can help to steer the simulation such that a
given error rate can be obtained with the smallest effort.

Fully adaptive simulations, possibly even using dynamic meshes that 
change from time step to time step, call for an enormous effort in
terms of implementation that is usually only given in academic
software codes (such as Gascoigne 3d,~\cite{Gascoigne}, which is used
in this 
work). Global climate models do not offer this flexibility. However, 
some models like
FESOM (\cite{Danilov2015}), MPAS (\cite{Ringler2013}) or ICON (\cite{Korn2017}) offer the possibility for
regional refinement in different zones. The error estimator can be
used for automatically selecting the proper refinement level of all
zones to reach the best accuracy on a discretization that is as coarse
as possible.

We start by describing a benchmark problem that has been introduced
by~\cite{MehlmannRichter2016newton}. While keeping the test case simple
(e.g. square domain) it features typical characteristics in terms of
the forcing and the parameters. Then, we present different numerical
studies on the error estimator. First we test its accuracy and
effectivity in terms of adaptive mesh control, then focusing on
possible cases for an integration of such techniques in climate models.

\subsection{Definition of a benchmark problem}\label{sec:bench}

For all test cases we consider the domain
  $\Omega=(0,500\unit{km})^2$.  At initial time $t=0$, the ice is at rest, $\vt_0=0$, the ice concentration
is constant $A=1.0$ and the ice height is a spatial variation around a
thickness of $H=\unit[0.3]{m}$.
\[
H^0(x,y)=\unit[0.3]{m}+\unit[0.005]{m}
\left(\cos\left(\frac{x}{\unit[25]{km}}\right)+ 
\cos\left(\frac{y}{\unit[50]{km}}\right)\right).
\]
A circular steady ocean current is described by
\[
\vt_\text{ocean}(x,y)= \unit[0.01]{m\cdot
  s^{-1}}\begin{pmatrix}y/\unit[250]{km}-1\\ 
  1-x/\unit[250]{km}
\end{pmatrix}. 
\]
The wind field mimics a cyclone and anticyclone that is diagonally
passing through the computational domain going back and forth
\[
v_\text{atm}(x,y,t) =\unit[15]{m\,s^{-1}}\omega(x,y)
R(\alpha)
\begin{pmatrix}
  x-m_x(t) \\ y-m_y(t)
\end{pmatrix},
\]
with the rotation matrix
\[
R(\alpha):=
\begin{pmatrix} 
  \cos(\alpha) & \sin(\alpha)\\
  -\sin(\alpha) & \cos(\alpha)
\end{pmatrix}.
\]
At initial time, the center of the cyclone is 
at $m_x(0)=m_y(0) = \unit[250]{km}$ and in 
$t\in [0,\unit[4]{days}]$ is moves to $m_x(4)=m_y(4) = \unit[450]{km}$
at constant speed. Then for $t\in [\unit[4]{days},\unit[12]{days}]$ it
goes back to $m_x(12)=m_y(12) = \unit[50]{km}$ where it turns again
towards $m_x(20)=m_y(20) = \unit[450]{km}$, and so on. When the
direction is towards the upper left, the convergence angle is set to $\alpha =
90^\circ-18^\circ$ and when it goes back towards the lower left we use 
$\alpha=90^\circ-9^\circ$. 
To reduce the wind
strength away from the center, we choose  $\omega(x,y)$ as
\[
\begin{aligned}
  \omega(x,y) &=
  \frac{1}{50}\exp\left(-\frac{r(x,y)}{\unit[100]{km}}\right),\\
  r(x,y) &= \sqrt{(x-m_x(t))^2+(y-m_y(t))^2}. 
\end{aligned}
\]
As functional of interest, we evaluate the average sea ice extent
within a subset $\Omega_2\subset\Omega$ of the domain
\begin{equation}\label{goal:seaice}
  \begin{aligned}
    \frac{1}{T}J_{A}(A)=\int_0^T\int_{\Omega_2}A(x,y,t) \,d(x,y)\,dt. 
  \end{aligned}
\end{equation}
We specify $\Omega_2$ for each test case. Similar measures are considered for 
sea ice model evaluations or model intercomparisons,
see~\cite{Stroeve2014} or \cite{Kwok2009}. The exact choice of the
subdomain $\Omega_2$ and also the time interval of interest $I=[0,T]$
will be specified in the different test cases.

\paragraph*{Solution of the nonlinear and linear systems}

The nonlinear problems resulting in each time step of the forward
simulation are solved with a modified Newton scheme that is described
in~\cite{MehlmannRichter2016newton}. The linear
systems within the Newton iteration and the linear problems of the dual system are solved with a GMRES method, preconditioned by a
geometric multigrid solver as
introduced in~\cite{MehlmannRichter2016mg}. The model is implemented  
in the software library Gascoigne 3d, see~\cite{Gascoigne}.


\subsection{Sharpness of the error estimate}\label{sec:num1}

In a first test case we evaluate the sharpness of the error
estimator, i.e. its capability of predicting an quantitatively exact
error value. As noted in the introduction to this section, this
scenario may not be of highest use in applications. However, it is an
important test case for the validation of the estimator itself.
For this first test case we use the short time interval
$I=[0,\unit[1]{day}]$ and the 
subdomain for measuring the functional 
$\Omega_2=(\unit[375]{km},\unit[500]{km})^2\subset \Omega$.
On a fine mesh with $h_{ref}=\unit[1]{km}$  and with the time step size
$k_{ref}=\unit[0.125]{h}$ we obtain the value
\begin{equation}\label{ref:seaiceA2}
  \tilde J_{A}:= J_A(U_{k_{ref},h_{ref}})=1.49907 \pm 10^{-5},
\end{equation}
in reference units, which corresponds to the average sea ice extent of 15615$\text{km}^2$.
We will take $\tilde J_A$ as reference value for the following computations.

\begin{table*}[t]
\begin{center}

  \begin{minipage}[t]{0.68\textwidth}
    \resizebox{\textwidth}{!}{
  \begin{tabular}{r r|c|c|c|ccc}
    \toprule
    \multicolumn{1}{c}{$h$}&    \multicolumn{1}{c|}{$k$}
    &$J_{A}(U^s_{k,h})$&$\tilde J_{A}-J_{A}(U^s_{k,h})$&$\eta_{k,h}$&$\eta_h$&$\eta_k$&$\eta_\beta$\\
\midrule
64 km&8 h&1.49763&1.44$\cdot 10^{-3}$&2.01$\cdot 10^{-3}$&1.20$\cdot 10^{-3}$&2.65$\cdot 10^{-3}$&1.58$\cdot 10^{-4}$\\
32 km&8 h&1.49788&1.19$\cdot 10^{-3}$&1.38$\cdot 10^{-3}$&1.21$\cdot 10^{-4}$&2.19$\cdot 10^{-3}$&4.40$\cdot 10^{-4}$\\
16 km&8 h&1.49797&1.10$\cdot 10^{-3}$&1.30$\cdot 10^{-3}$&6.72$\cdot 10^{-5}$&2.10$\cdot 10^{-3}$&4.38$\cdot 10^{-4}$\\
8 km &8 h&1.49802&1.05$\cdot 10^{-3}$&1.25$\cdot 10^{-3}$&4.11$\cdot 10^{-5}$&2.03$\cdot 10^{-3}$&4.21$\cdot 10^{-4}$\\
\midrule
64 km&4 h&1.49833&7.43$\cdot 10^{-4}$&9.52$\cdot 10^{-4}$&6.28$\cdot 10^{-4}$&1.21$\cdot 10^{-3}$&6.12$\cdot 10^{-5}$\\
32 km&4 h&1.49849&5.80$\cdot 10^{-4}$&6.16$\cdot 10^{-4}$&8.53$\cdot 10^{-5}$&1.02$\cdot 10^{-3}$&1.25$\cdot 10^{-4}$\\
16 km&4 h&1.49856&5.15$\cdot 10^{-4}$&5.72$\cdot 10^{-4}$&4.41$\cdot 10^{-5}$&9.70$\cdot 10^{-4}$&1.30$\cdot 10^{-4}$\\
8 km &4 h&1.49858&4.87$\cdot 10^{-4}$&5.51$\cdot 10^{-4}$&2.47$\cdot 10^{-5}$&9.44$\cdot 10^{-4}$&1.32$\cdot 10^{-4}$\\
\midrule
64 km&2 h&1.49863&4.39$\cdot 10^{-4}$&5.67$\cdot 10^{-4}$&5.48$\cdot 10^{-4}$&5.67$\cdot 10^{-4}$&2.04$\cdot 10^{-5}$\\
32 km&2 h&1.49876&3.10$\cdot 10^{-4}$&2.94$\cdot 10^{-4}$&7.01$\cdot 10^{-5}$&4.82$\cdot 10^{-4}$&3.67$\cdot 10^{-5}$\\
16 km&2 h&1.49881&2.59$\cdot 10^{-4}$&2.67$\cdot 10^{-4}$&3.58$\cdot 10^{-5}$&4.59$\cdot 10^{-4}$&3.95$\cdot 10^{-5}$\\
8 km &2 h&1.49883&2.37$\cdot 10^{-4}$&2.54$\cdot 10^{-4}$&1.93$\cdot 10^{-5}$&4.47$\cdot 10^{-4}$&4.19$\cdot 10^{-5}$\\
\bottomrule
  \end{tabular} 
    }
  \end{minipage}
  \begin{minipage}[t]{0.3\textwidth}
    \includegraphics[width=0.98\textwidth]{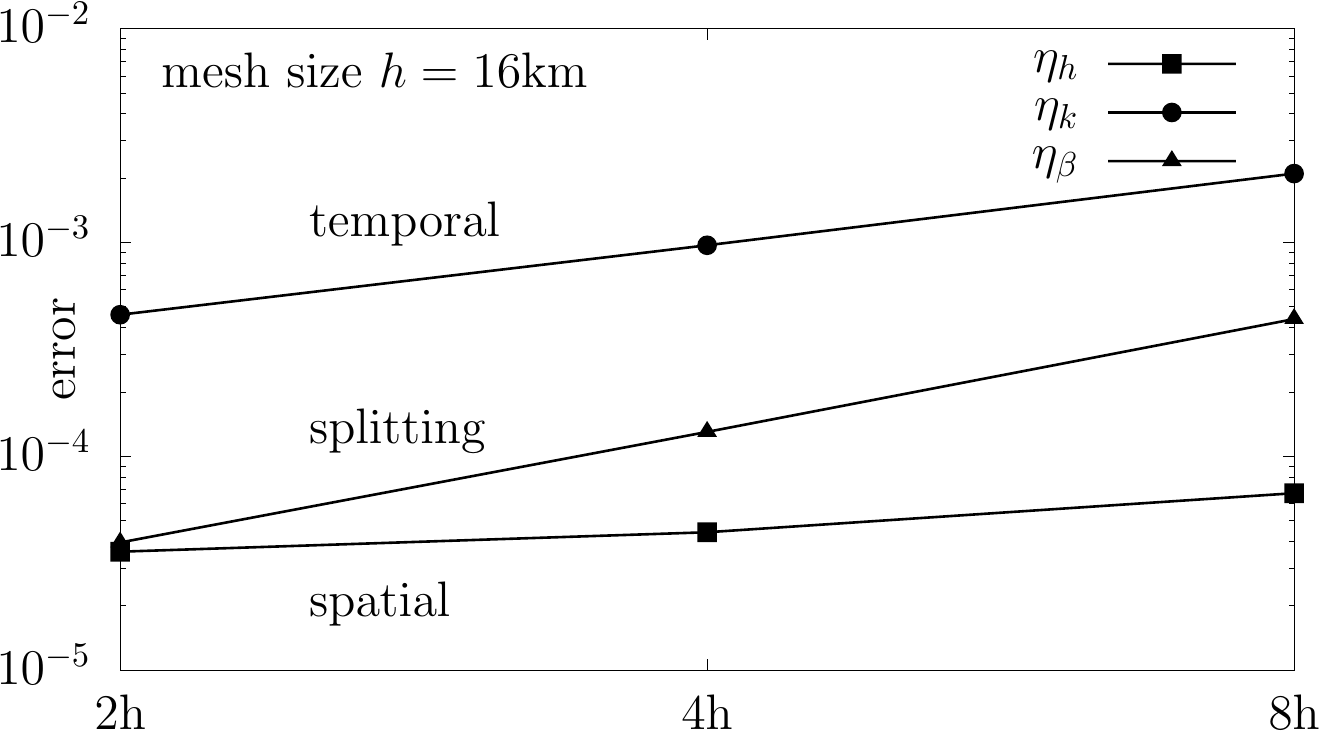}
      
    \includegraphics[width=\textwidth]{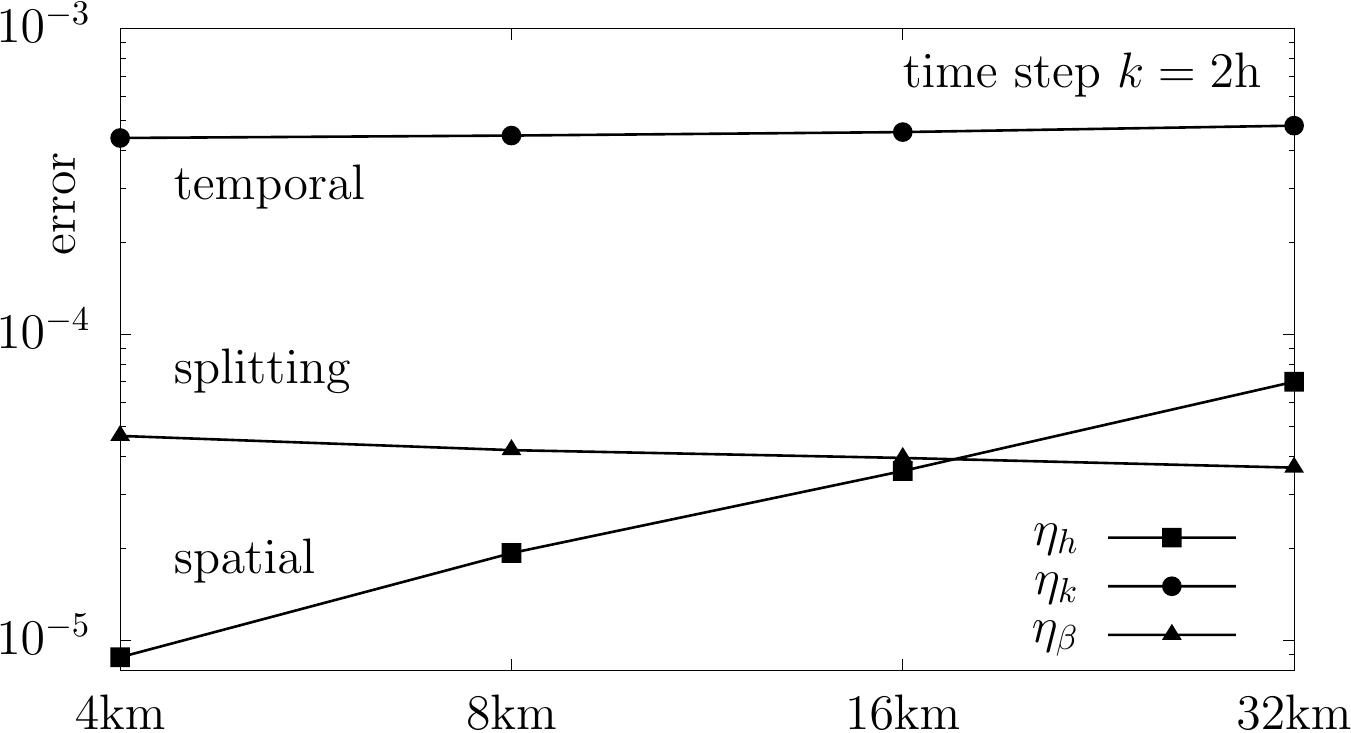}
  \end{minipage}
  \caption{On a sequence of spatially and temporally refined meshes we
    compare the functional output $J_A(U_{k,h}^s)$ to the reference
    value $\tilde J_A\approx 1.49907$ and indicate the error $\tilde
    J_A-J_A(U_{k,h}^s)$, the total error estimator
    $\eta_{k,h}=\frac{1}{2}(\eta_h+\eta_k+\eta_{\beta})$ and
    its contributions attributing the spatial discretization error
    $\eta_h$, 
    the temporal discretization error $\eta_k$ and the splitting error
    $\eta_\beta$. We observe that the complete error estimate
    $\eta_{k,h}$ is very close to the real error and that the temporal
    error is dominating on fine spatial meshes. On the right, we show
    the composition of the error estimator into spatial error
    $\eta_h$, temporal error $\eta_k$ and splitting error
    $\eta_\beta$ for a fixed spatial mesh and varying time step sizes
    (top) and for a fixed temporal mesh and varying spatial mesh sizes
    (bottom). }
  \label{tab:func}
\end{center}
\end{table*}

\begin{figure}[t]
  \begin{center}
    \includegraphics[width=0.5\textwidth]{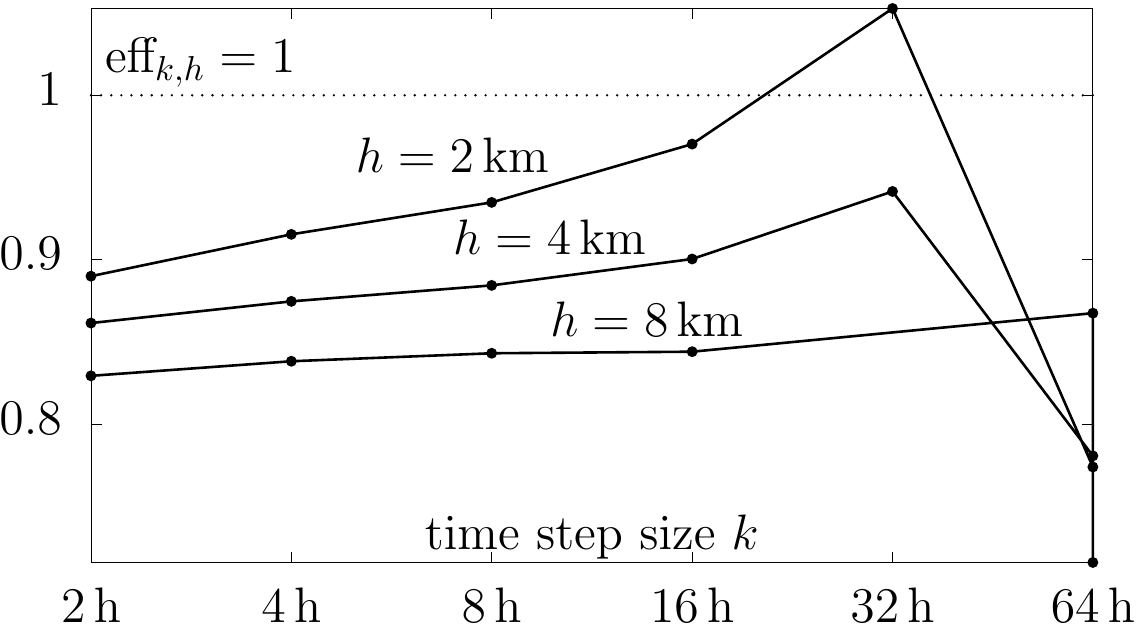}
  \end{center}
  \caption{\label{fig:eff_seaice}
    Effectivity index $\text{eff}_{k,h}$ of the error estimator for
    different temporal and spatial discretizations. The estimator is
    highly accurate with only about $20\%$ of overestimation.}
\end{figure}

In Table \ref{tab:func}, we evaluate the functional error $|\tilde
J_{A} -J_{A}(U^s_{k,h})|$ and the error estimator given
by~(\ref{error_est}), which we denote by $\eta_{k,h}$.
In the plots on the right side of Table~\ref{tab:func}, we show the composition of the error estimator into
  spatial error, temporal error and splitting error. The upper figure
  shows that, for a fixed spatial mesh with $h=16$km, linear
  convergence is obtained for $\eta_k$ and $\eta_\beta$ under temporal
  refinement, whereas the spatial error naturally stagnates. The lower
  figure shows in a similar fashion that, for a fixed time mesh with
  $k=2$h, spatial refinement results in linear convergence in space,
  whereas temporal error and splitting error do not get smaller. This
  clear decomposition of the error estimator into spatial and temporal
  contributions will allow us to design an adaptive algorithm that
  efficiently controls the discretization in order to balance temporal
  and spatial error, see Section~\ref{sec:balance}. Overall, this
  example shows a dominance of the temporal error, which is due to the
short simulation time of $1$ day, where nearly no kinematic features
appear. To validate the
accuracy of the error estimator we introduce the \emph{effectivity index}
\begin{equation}\label{eff_inex}
  \text{eff}_{k,h}:=\frac{\tilde J_A-J(U^s_{k,h})}{\eta_{k,h}},
\end{equation}
which measures the sharpness of the estimate. If this index is
  close to one, true error and estimator are very close
  $\eta_{k,h}\approx \tilde J_A-J(U^s_{k,h})$ such that the estimate
  is very accurate. If the index is much larger than one, the
  estimator overestimates the true error, if it is much smaller than
  one, the estimator underestimates the error.
  In Figure \ref{fig:eff_seaice}, we plot this effectivity index and
  find $0.75\le \text{eff}_{k,h} \le 1.1$ which indicates that
  the error estimator is  highly accurate, in particular for
  increasing spatial  and temporal resolutions. 

The last three columns of Table~\ref{tab:func} show the decomposition of
the error estimator into spatial, temporal and splitting part as
described in Section~\ref{sec:decomp}. These 
values show a dominance of the temporal error over the spatial error
and to lesser degree also over the splitting error. In space, the
estimator values $\eta_h$ also clearly demonstrate linear
convergence in $h$, which is expected for linear finite elements. We
do not observe this convergence order in the overall error, as it is
dominated by the other two parts. The second test case in Section \ref{sec:adaptive} shows a dominance of the spatial 
error. The 
dominating \emph{temporal residual error}  stems from  the short
simulation time of $T=\unit[1]{day}$. Our findings coincide with the
analysis of~\cite{Lemieux2014} where the temporal error also 
dominates the splitting error in a one day simulation. A test case running for $\unit[33]{days}$ discussed in Section~\ref{sec:adaptive} shows a balanced distribution of spatial and temporal errors.


\subsection{Balancing error contributions}\label{sec:balance}

A simple application of the decomposition of the error estimator into
spatial error, temporal error and splitting error is to balance
the different error contributions by the following  algorithm:
\begin{algo}[Balancing errors]\label{algo:balance}
  Given an initial time step size $k$ and mesh size $h$. Iterate:
  \begin{enumerate}
  \item Solve the sea ice problem $u_{k,h}\in V_{k,h}$
  \item Estimate the error according to Algorithm~\ref{algo:feedback}
  \item Split the error estimate
    $\eta_{k,h}:=\frac{1}{2}(\eta_h+\eta_k+\eta_\beta)$ 
  \item If $\eta_k+\eta_\beta>2 \eta_h$ refine time step $k\mapsto
    \frac{k}{2}$\\
    If $\eta_h>2(\eta_k+\eta_\beta)$ refine spatial mesh $h\mapsto
    \frac{h}{2}$\\
    Otherwise refine in space and time $k,h\mapsto \frac{k}{2},\frac{h}{2}$\\
  \end{enumerate}
\end{algo}
Here, we have attributed the splitting error $\eta_\beta$ to the
temporal error. We refine only spatially (or temporally) if this error
contribution is twice as large as the other part. If the errors are
already close to each other, we refine in space and in time. This
strategy can be extended to include further error contributions. In
coupled ice-ocean simulations one could balance the errors of the
ocean component and the ice component.

We virtually perform a simulation based on
Algorithm~\ref{algo:balance} by processing the results from
Table~\ref{tab:func}. Starting with $h=\unit[64]{km}$ and
$k=\unit[8]{h}$ it holds $\eta_k+\eta_\beta=2.71\cdot
10^{-3}>2\eta_h=2.40\cdot 10^{-3}$ (compare the first line of
Table~\ref{tab:func}). Hence, we refine in time only and proceed with
$h=\unit[64]{km}$ and $k=\unit[4]{h}$. Again, it holds
$\eta_k+\eta_\beta=1.27\cdot 10^{-3}>2\eta_h=1.26\cdot 10^{-3}$ such
that we once more refine in time only, resulting in $h=\unit[64]{km}$
and $k=\unit[2]{h}$. This third simulation yields $\eta_k+\eta_\beta =
5.87\cdot 10^{-4}\approx \eta_h=5.48\cdot 10^{-4}$ and we would
continue by refining  both  in time and space.

\begin{figure*}[t]
  \begin{tabular}{c c c}
    \includegraphics[height=0.3\textwidth]{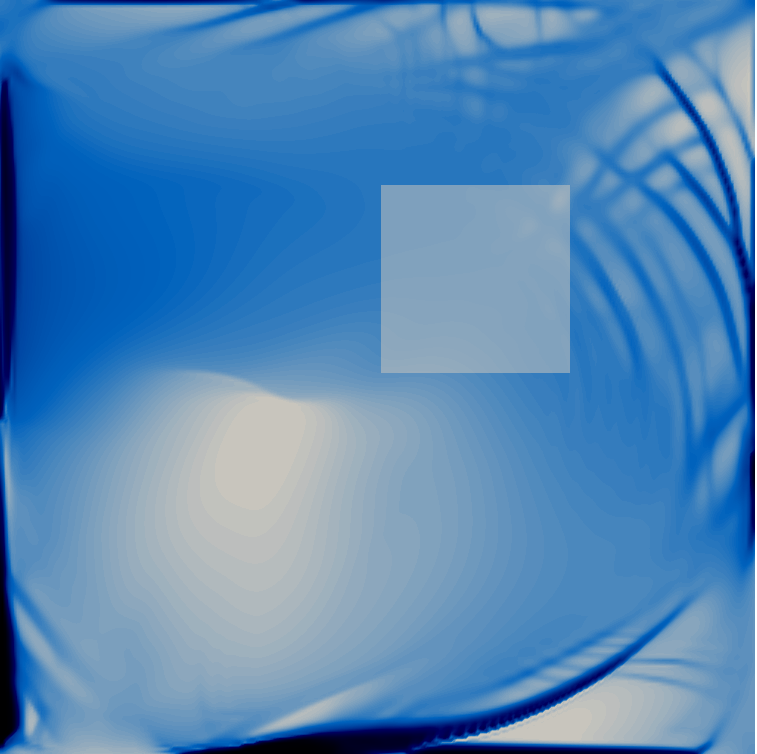}&
    \includegraphics[height=0.3\textwidth]{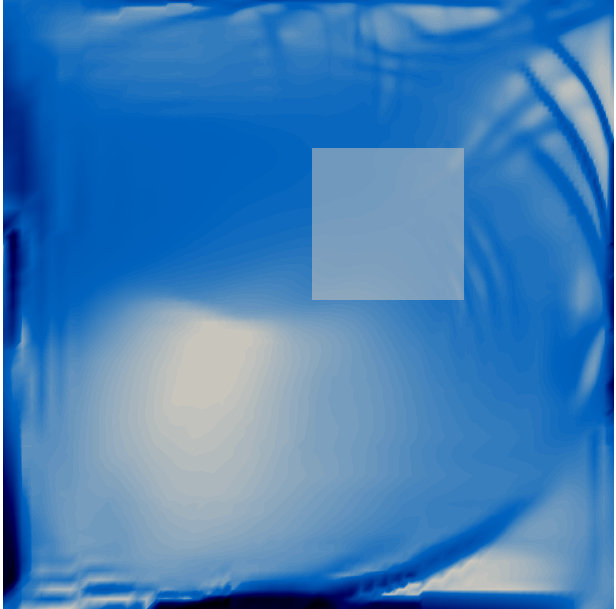}&
    \includegraphics[height=0.3\textwidth]{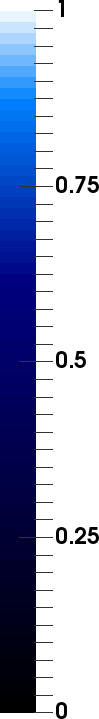}
  \end{tabular}
  \caption{Solution at day 15. On the left a uniform mesh with
    horizontal spacing $\unit[2]{km}$ and on the right an
    adaptive mesh using non-uniform mesh spacing between
    $\unit[64]{km}$ up to $\unit[2]{km}$ is used. The highlighted 
    area is the domain $\Omega_2$ where we measure the average sea ice
    extent.}
  \label{fig:ada1}
\end{figure*}

The natural alternative to this procedure would be a uniform
refinement in space and in time whenever the accuracy is not
sufficient. To compare the complexity of both approaches we assume
that the algorithm scales optimally, i.e. linear in the number of time
steps ${\cal O}(k^{-1})$ and linear in the number of mesh elements
given by ${\cal O}(h^{-2})$.\footnote{Linear complexity w.r.t. spatial
  refinement is in principal possible by using multigrid methods for
  the solution of the linear systems,
  see~\cite{MehlmannRichter2016mg}. Due to the increasing impact of
  the nonlinearity on highly resolved simulations, the assumption of
  linearity turns out to be too optimistic. The savings from
  adaptivity by using smaller meshes would even be more drastic if a
  realistic estimate of the effort would be available.}%
Altogether we use the
simple model $E(k,h) = C k^{-1}h^{-2}$ to measure the effort of one
simulation. For simplicity, the constant is set to $C=64^2\cdot 8$. 
Three steps of uniform refinement result in the 
effort
\[
E(8,64)+E(4,32)+E(2,16)=
32\,768\Big(\frac{1}{64^2\cdot 8}+
\frac{1}{32^2\cdot 4}+
\frac{1}{16^2\cdot 2}\Big)=73
\]
whereas the balancing algorithm yields
\[
  E(8,64)+E(4,64)+E(2,64)=
  32\,768\Big(\frac{1}{64^2\cdot 8}+
  \frac{1}{64^2\cdot 4}+
  \frac{1}{64^2\cdot 2}\Big)=7,
\]
which is only $10\%$ of the effort for the uniform standard
approach. On the final mesh, the balancing algorithm yields the error
$4.39\cdot 10^{-4}$ compared to $2.59\cdot 10^{-4}$ that would be
obtained by using uniform refinement in space and time (at 10 times
the cost).


\subsection{Adaptive mesh control and steering of regional
  refinement}\label{sec:adaptive}

   
In a second test case we consider a longer time horizon of
$T=\unit[33]{days}$ and an initial sea ice concentration of $A=0.9$. The wind field described in
Section~\ref{sec:bench} passes the domain several times and typical
kinematic features appear. Figure~\ref{fig:ada1} shows the
sea ice ice concentration at day 15. The left figure gives a result on
a uniform discretization, whereas the right plot belongs to the
corresponding result on a locally refined mesh with about 3 times less
  unknowns. The highlighted quadrilateral area is the domain
$\Omega_2=(\unit[250]{km},\unit[375]{km})^2$ where we evaluate the
average sea ice extent
\begin{equation}\label{f2}
  J(A) = \frac{1}{\unit[33]{days}} \int_0^{\unit[33]{days}}
  \int_{\Omega_2}A(x,t)\,\text{d}{x}\,\text{d}t.
\end{equation}

\begin{figure*}[t]
  \begin{center}
    \includegraphics[width=\textwidth]{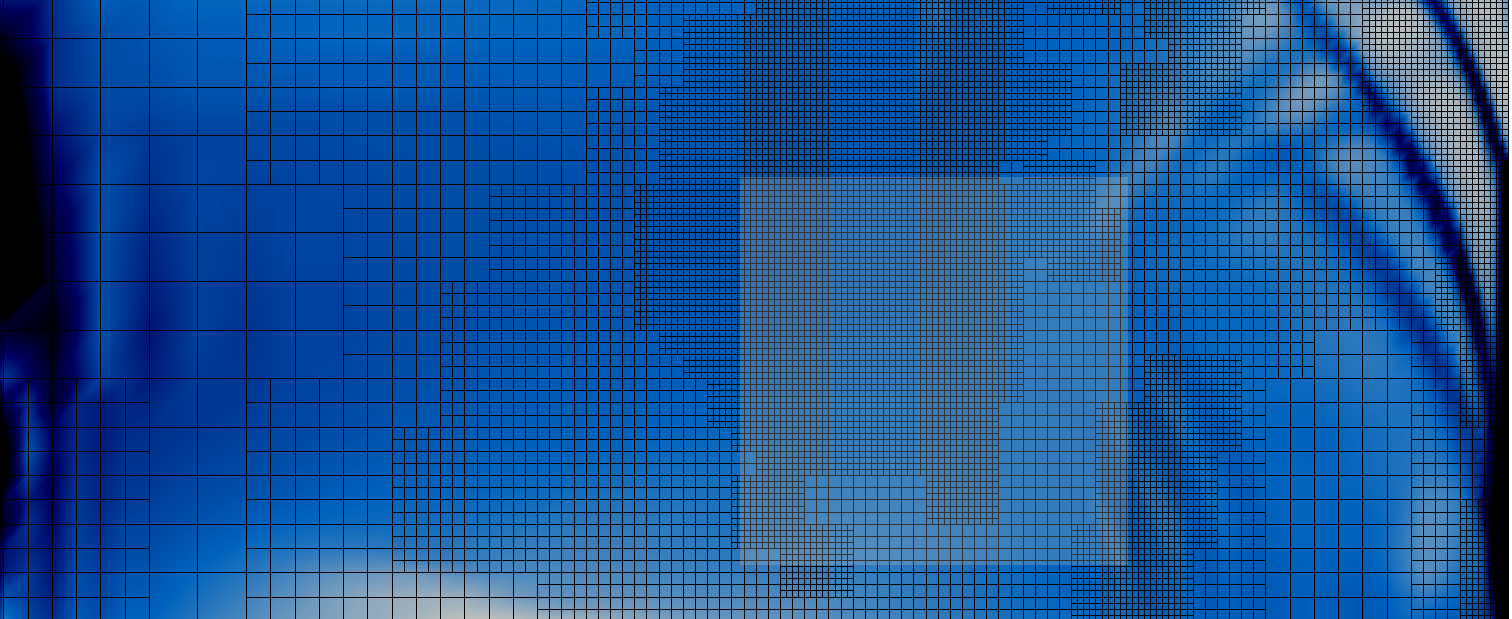}
  \end{center}
  \caption{Zoom into the adaptive mesh in the area of interest.}
  \label{fig:ada2}
\end{figure*}

We focus on the spatial convergence and use the step size
$k=\unit[0.5]{h}$ for all tests. We consider the contribution
of the spatial discretization error $\eta_h$ only and use
Algorithm~\ref{algo:feedback} for identifying optimal finite element
meshes to yield small errors on meshes that are as coarse as
possible. Two different strategies are investigated. First, we use a
fully local adaptive mesh concept, where mesh elements
$K\in\Omega_h$ are refined into four smaller quadrilaterals, if the
local error contribution $\eta_K$ is larger than the average error,
i.e.
\begin{equation}\label{refine}
  \eta_K>\gamma\cdot \bar \eta\quad\Rightarrow\quad\text{refine }K,\quad
  \bar\eta\coloneqq  \frac{1}{|\Omega_h|} \sum_{K'\in\Omega_h} \eta_{K'}
\end{equation}
where $\gamma\approx 1$ is a constant to fine-tune the refinement
procedure. Usually we take $\gamma=2$. 
The results are shown in Figure~\ref{fig:ada2}. We focus on a small
region around the area of interest $\Omega_2$ to better highlight the
mesh that has been generated by the error estimator.
It is not necessary to resolve the complete region of
interest $\Omega_2$. Instead, parts outside of this region also
have to be resolved to get the correct transport of
information. Further, the resulting meshes appear rather tattered and
non-symmetric. It is a typical feature of the error estimator that the
resulting meshes are rather non-intuitive. We refer also
to~\cite{BeckerRannacher2001} with several examples showing that
meshes obtained by a posteriori error estimators are superior to
manually adjusted refinements. 
The solution on the adaptive discretization still shows
kinematic features. These however are less distinct in comparison to
the global discretization, compare Figure~\ref{fig:ada1}, where we show
the uniform result (left) and the adaptive one (right) side by
side. We have to keep in mind that the goal of our estimator is not to
detect features but to predict the average sea ice
extent~(\ref{f2}). The dual weighted residual estimator consists of residuals that measure the exactness of the solution and of the adjoint weights, which measure the sensitivity with respect to the goal functional. We only refine, if both the local residual and the local sensitivity information indicate a large and relevant error. 
 Adaptation of the mesh is guided by finding the best mesh allocation, that minimize the approximation error in the evaluation of the functional of interest $J(\cdot)$. 
 Hence, typical structures of the solution are not good indicators of the quality of the adaptative discretization.



\begin{figure}[t]
  \begin{minipage}{0.5\textwidth}
  \includegraphics[width=\textwidth]{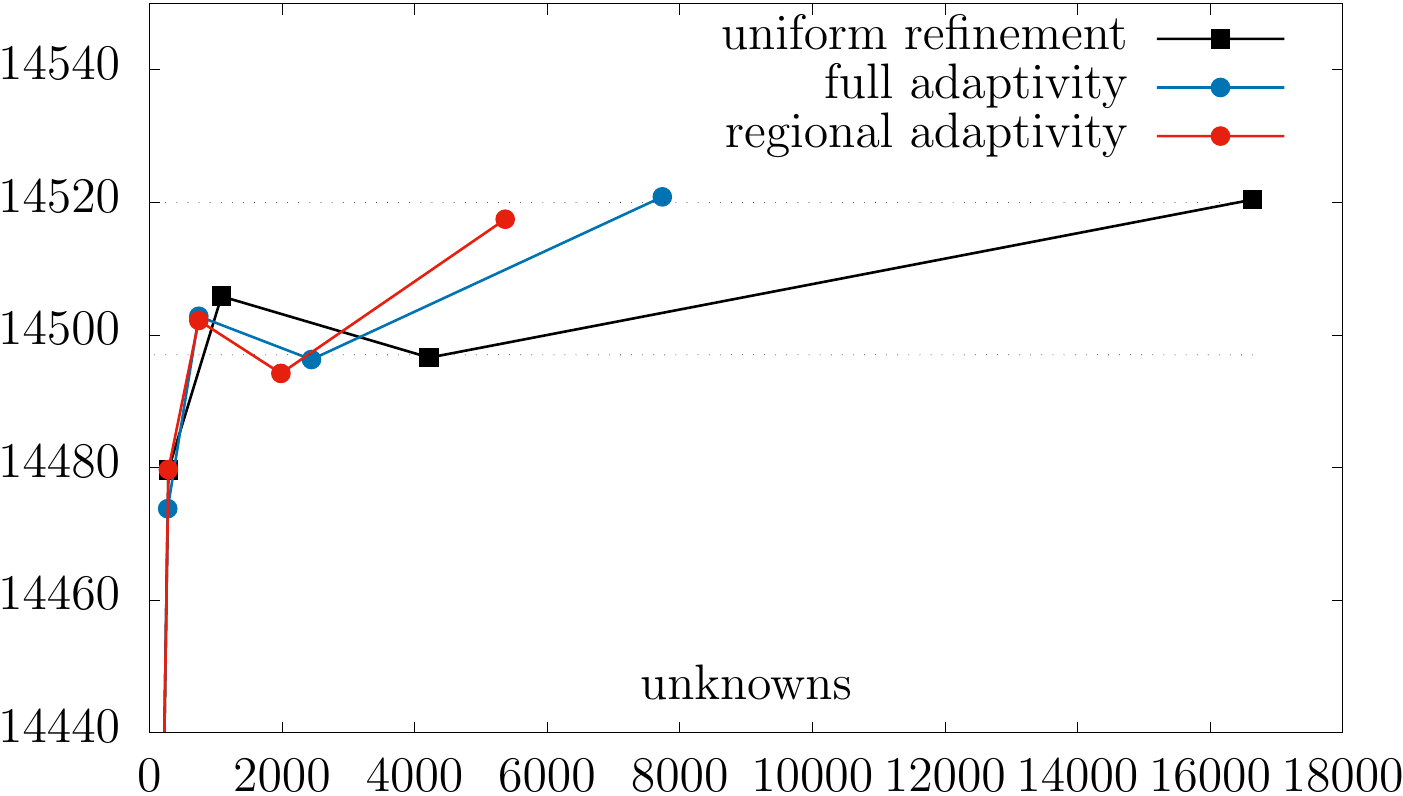}
  \end{minipage}\hspace{0.05\textwidth}
  \begin{minipage}{0.44\textwidth}
    We show (black line with squares) the
    results for global refinement, (blue line with bullets) results
    for fully adaptive meshes and (red line with diamonds) the results
    for adaptive meshes based on regions. For better comparison of the
    different approaches we add dotted lines indicating the functional
    levels obtained on uniform meshes.
  \end{minipage}
  \caption{Spatial discretization error vs. the number of unknowns
    using adaptive meshes. }
  \label{fig:ada}
\end{figure}

Instead,  we show in Figure~\ref{fig:ada} the resulting average sea ice
extent on uniform and adaptive discretizations. We observe that the
adaptive algorithm is able to capture exactly the same dynamics as the
uniform discretization, but, by using fewer unknowns and hence on a
significantly reduced problem size. The benefit of adaptivity is the
omitting of unnecessary refinements. Due to the nonlinearity of the
viscous-plastic sea ice model and the appearance of features in the
solution that numerically nearly resemble discontinuities, we do not
observe a monotonic convergence of the functional $J(A_h)$ for $h\to
0$. Therefore, Figure~\ref{fig:ada} shows the functional output
itself. With $|\Omega_2| = \unit[15\,625]{km^2}$ these numbers
correspond to an ice cover of approximately $93\%$. 

As adaptive meshes are able to give similar quality in the goal
functional on smaller meshes, the computational efficiency is
significantly reduced. The two simulations shown in
Figure~\ref{fig:ada1} both belong to discretizations with a minimum
mesh size of 2 km. The overall computational time for solving the 33
day-test case on uniform meshes was 81 hours (about 3 days). The
corresponding simulation on the adaptive mesh took 13 hours. Adding 
the complete overhead of the error estimator (computation of the
adjoint problem and evaluation of the residuals), the computational
time sums up to 20 hours, four times less than the fully uniform
simulation.\footnote{All computations have been carried out on a laptop using the single core performance of a Core i5-6360U CPU at 2.0 GHz.}

As discussed before, usual large scale climate models do not allow for
fully adaptive meshes that call for a large technical overhead in
terms of implementation, in particular when it comes to efficient
realizations on parallel computers. However, several models allow for
selecting local regions of higher resolution. These are usually
hand-picked. Here we discuss a second possible use of the a posteriori error
estimator for optimally tuning the mesh sizes in predefined local
regions. We split the domain $\Omega=(0,\unit[500]{km})^2$ into
$16=4\times 4$ uniform local regions. Then, we proceed similar
to~(\ref{refine}) but first sum all error indicators $\eta_{K;h}$ that
belong to each of the 16 regions. Refinement is not carried out
element-by-element, but for the complete region that has been
selected. The corresponding results are shown in the red line of
Figure~\ref{fig:ada}. 

The course of the functional values obtained with regional refinement
is similar to the uniform and the fully adaptive case. Even
slightly less unknowns are chosen as compared to the fully adaptive
case. However the functional values are a bit off and the regional
refinement method is not able to completely match the uniform
discretization. The reason is found in the averaging of the error
estimators to the 16 regional blocks. Only if the average is above a
certain limit, the complete block is refined. By reducing the
parameter $\gamma$ in~(\ref{refine}), more refinement could be
achieved. Too low values of $\gamma$ might however result in
unnecessary overrefinement.

\begin{figure}[t]
  \begin{center}
    \includegraphics[width=\textwidth]{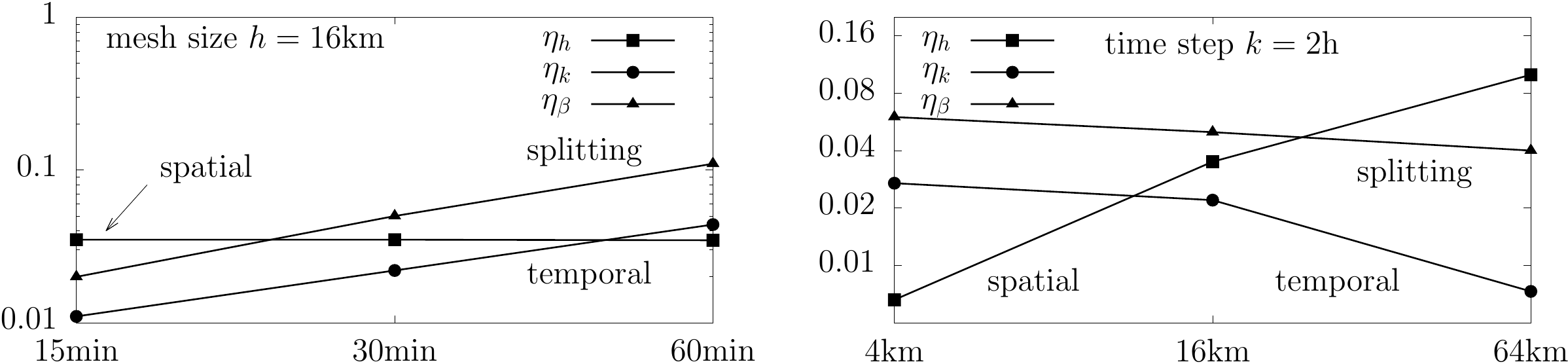}
  \end{center}
  \caption{Decomposition of the error estimator into spatial, temporal
    and splitting error. Left: mesh size $h$ fixed and right: time
    step size $k$ fixed. }
  \label{fig:decomp-long}
\end{figure}

Finally, we show in Figure~\ref{fig:decomp-long} the decomposition of
the error estimator into spatial error $\eta_h$, temporal error
$\eta_k$ and splitting error $\eta_\beta$. The corresponding study for
the short term test case is presented on the right of
Table~\ref{tab:func}. While the temporal error was dominating there,
the long time example shows a more prominent spatial and splitting error. 


 
\section{Discussion and Conclusion}\label{sec:con}
In this paper we introduced the first error estimator for the standard
model describing the sea ice dynamics. The error estimator is derived
for a general class of coupled  non-stationary partial differential
equations that are solved with a {partitioned solution approach}. It
is based on the concept of the {dual weighted residual method} that
has been introduced by~\cite{BeckerRannacher2001}. The error estimator
consists mainly of two parts, the \emph{primal} and \emph{dual residual error}
that arise in the framework of the dual weighted residual method, and
for the first time, 
an additional \emph{splitting error} which stems from the application of
the partitioned solution approach, is considered. In order to derive the error
estimator for the sea ice model, we reinterpret the usual implicit
Euler formulation as a variational space-time Galerkin approach. 

We numerically evaluated this new error estimator on an idealized test
case and measured the sea ice extent in a subdomain of interest. The
temporal discretization error dominates the 
overall numerical error on all considered mesh resolutions in the one day simulation. This might
be due to the short simulation time and it coincides with the findings
of~\cite{Lemieux2014}. Considering a 33 day simulation  on coarse meshes the spatial error is dominant. With increasing mesh resolution the splitting error becomes the most important error source.

The error estimator is highly accurate as we observe an efficiency
index close to 1. Despite the very strong nonlinearity of the sea ice
model this means that the DWR estimator is a useful measure in sea ice
simulations. 

We discussed several approaches
how this error estimator can be used to speedup the sea ice component
in global climate models. First, the error estimator can be applied
for a balancing of different error contributions, namely the spatial
and the temporal discretization error as well as the error that comes
form partitioning the coupled system. This approach can be extended to include further
fields, like a coupled ocean-ice simulation. Second, we demonstrate
how the error estimator can be used to control the mesh size of models
that allow for a regional sampling at higher resolution. An automatic
feedback approach guides the simulation to an optimally balanced mesh
and allows for significant savings in terms of computational time.

Based on the work of~\cite{BraackErn2003} one could extent the error estimator to also include a model error. One promising application is to consider the adaptive EVP model (see~\cite{Losch2016}) as an approximation to the VP model. The discrepancy between VP and EVP model can be included in terms of residual evaluations such that a balancing of discretization error and model error will result in an effective stopping criteria for the adaptive EVP iteration. Similarly, iteration errors coming from approximate Picard iterations can be taken care of by assuming a further disturbance of the Galerkin orthogonality. Details are discussed in~\cite{MeidnerRannacherVihharev2009}.

The main technical difficulty for realizing the error estimator is the
implementation of the dual problem,  that runs backward in
time and that has a reversed partitioning structure. Such adjoint
  solutions are also essential in variational data assimilation and
  some climate models offer implementations.
The concept of the dual weighted residual estimator is very
flexible, with the main prerequisite of casting the problem and
discretization into a variational Galerkin formulation. We have
considered one typical error functional measuring the average ice
extent, but further error measures are easily realized.

\medskip
\noindent\textbf{Acknowledgment.}
The work of Carolin Mehlmann has been supported by the \emph{Deutsche
  Bundesstiftung Umwelt}. The work of both authors is funded by the
Deutsche Forschungsgemeinschaft (DFG, German Research Foundation) -
314838170, GRK 2297 MathCoRe. We thank the anonymous reviewers for their effort that helped to improve the quality of this manuscript.


\appendix
\section{Appendix}
\subsection{Temporal Galerkin Discretizations}\label{app:gd}

\begin{figure}[t]
\begin{center}
 \includegraphics[width=0.4\textwidth]{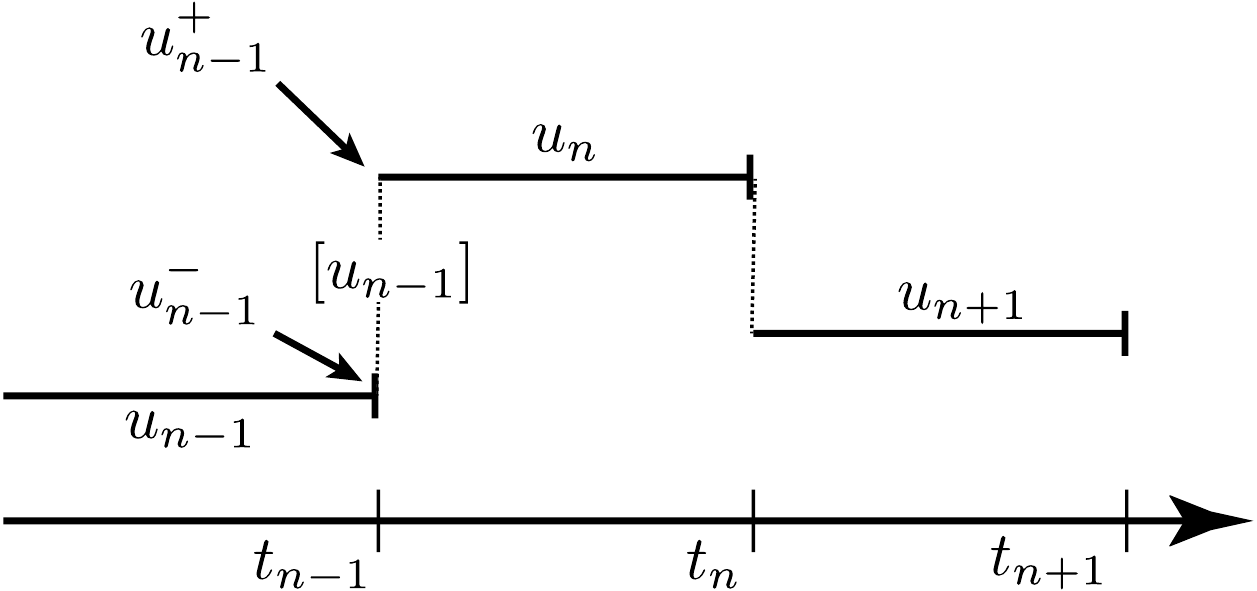}
 \end{center}
 \caption{\label{Def:spruenge}Visualization of the temporal jump of piecewise constant functions $u_k$ with $u_n=u_k\big|_{(t_{n-1},t_n]}$ at time point $t=t_n$. } 
\end{figure}

Since variational formulations in space and time are the basis for the dual weighted residual estimator we briefly describe the relation between the classical backward Euler time stepping method and the temporal dG(0) discretization used for estimating the error. Considering the ode $u'(t)=f\big(t,u(t)\big)$ with $u(0)=t_0$, its variational formulation on a partitioning $0=t_0<t_1<\dots<t_N=T$ is given by
\begin{equation}\label{A1}
  A(u,\phi)\coloneqq 
  \sum_{n=1}^N \Big\{\int_{t_{n-1}}^{t_n} \Big(u'(t)-f\big(t,u(t)\big)\Big)\cdot \phi(t)\,\text{d}t + [u(t_{n-1})]\cdot \phi(t_{n-1})^+ \Big\}= 0,
\end{equation}
where $[\cdot]$ denotes the jump of the possibly discontinuous function $u(t)$, compare Section~\ref{sec:var}. Each smooth solution $u'(t)=f(t,u(t))$ naturally satisfies this variational formulation $A(u,\phi)=0$. If we discretize~(\ref{A1}) with piecewise constant functions in $u$ and $\phi$, i.e. $u_n=u\big|_{(t_{n-1},t_n]}\in\mathds{R}$ and $\phi_n=\phi\big|_{(t_{n-1},t_n]}\in\mathds{R}$ the sum in~(\ref{A1}) decouples into discrete time steps and the integral in time can be computed exactly with the box rule
\begin{equation}\label{A2}
  \Big(-(t_n-t_{n-1})f(t_n,u_n)+ (u_n-u_{n-1})\Big)\phi_n = 0,
\end{equation}
which, after dividing by $t_n-t_{n-1}$ and by using $\phi_n=1 \in\mathds{R}$, gives the backward Euler method. Hereby we can state, that the dG(0) Galerkin discretization $A(u,\phi)=0$ defined in~(\ref{A1}) and the backward Euler method~(\ref{A2}) are equivalent,\footnote{Since the derivation of~(\ref{A2}) from~(\ref{A1}) relies on the evaluation of the temporal integrals with the box rule, equivalence only holds, if all integrals are computed exactly. This is the case for autonomous equation, where $f=f(u(t))$ does not explicitly depends on $t$, but not in the general case.} hence, $A(u_k,\phi_k)=0$ also holds for the backward Euler solution $u_n$, interpreted as piecewise linear function.  
Figure~\ref{Def:spruenge} shows the discrete Galerkin solution as a piecewise constant function with discontinuities at the time steps $t_n$.  

For the formulation of the error estimator it is nevertheless essential to have the variational formulation~(\ref{A1}) in mind, since the estimator relies on \emph{weighted residuals}, i.e. on evaluations of~(\ref{A1}) with test functions that come from a space of higher degree (e.g. $\phi_h=i_{k}^{(1)}z_h - z_h$, polynomials of degree 1). For higher order degree polynomials the box rule is not longer exact for evaluating the integrals and we must indeed use the variational formulation.

\subsection{Galerkin orthogonality and partitioned solution}\label{app:go}

\begin{figure}[t]
  \begin{center}
    \includegraphics[width=0.4\textwidth]{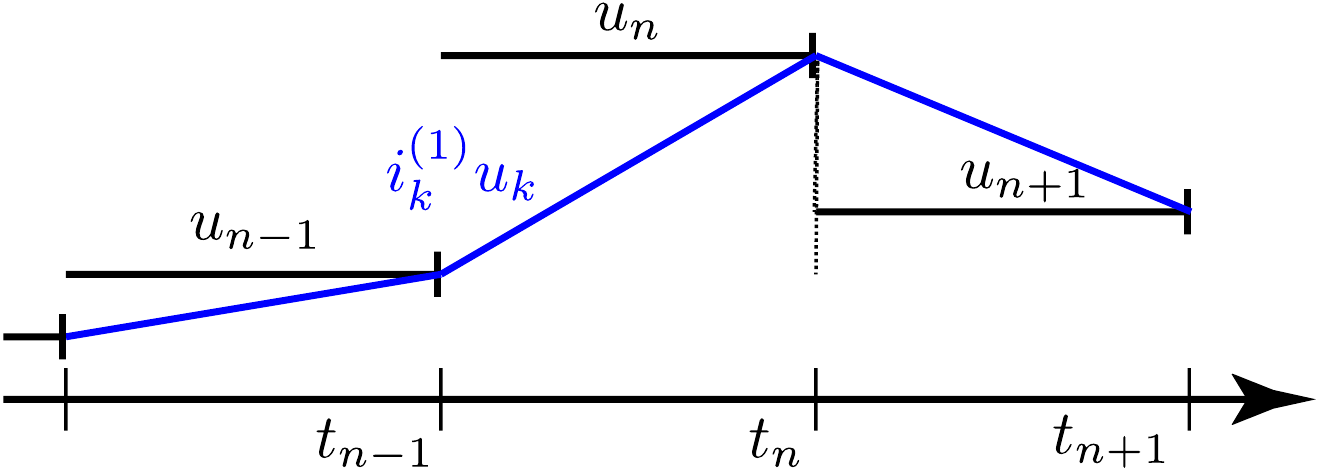}
    \caption{Reconstruction of the piecewise constant $u_k$ by piecewise linear polynomials for the evaluation of the weights.}
    \label{fig:recon}
  \end{center}
\end{figure}

Both the real solution $u(t)$ to $u'(t)=f(t,u(t))$ and the backward Euler approximation $u_k$ satisfy the variational formulation, $A(u,\phi)=0$ and $A(u_k,\phi_k)=0$, respectively. The difference is the choice of test functions $\phi$. While $u$ satisfies the variational problem for \emph{all} test functions, $A(u_k,\phi_k)=0$ holds only for piecewise constant functions $\phi_k$. For these, Galerkin orthogonality is satisfied
\begin{equation}\label{A3}
  A(u,\phi_k)-A(u_k,\phi_k)=0.
\end{equation}
Galerkin orthogonality plays an important role in the derivation of the dual weighted residual estimator, which, for simple linear problem takes the form
\begin{equation}\label{A4}
J(u)-J(u_k) = A(u-u_k,z),
\end{equation}
where $z$ is the adjoint solution to $A(\psi,z)=J(\psi)$ for all test functions $\psi$. One problem, discussed in Remark~\ref{rem:weights} is the presence of the unknown solutions $u$ and $z$ that serve as \emph{weights} of the estimator. By Galerkin orthogonality~(\ref{A3}) we can add to~(\ref{A4}) an interpolation of the adjoint solution $i_k z$
\begin{equation}\label{A5}
J(u)-J(u_k) = A(u-u_k,z-i_k z) = -A(u_k,z-i_k z).
\end{equation}
These new weights $z-i_k z$ are still not known (since they involve $z$). However, they can be efficiently approximated by reconstructing the interpolation error $z-i_k z$ in a space of higher order (linear functions instead of constants), again, see Remark~\ref{rem:weights} and we obtain the approximation
\begin{equation}\label{A6}
J(u)-J(u_k) \approx -A(u_k,i_{k}^{(1)} z_k - z_k). 
\end{equation}
This is the most simple \emph{primal form} of the DWR method which is valid for linear models and linear functionals, compare~\cite{BeckerRannacher2001}. It only consists of the primal residual $A(u_k,\cdot)$ weighted with the dual interpolation error $i_{k}^{(1)}z_k - z_k$. In Section~\ref{sec:error_estimator} we must apply the more general form of the DWR method which also translates to nonlinear models like the sea ice momentum equation. Here, additional adjoint residuals, weighted with the primal interpolation errors $i_{k}^{(1)}u_k-u_k$ appear. In Figure~\ref{fig:recon} we illustrate this discrete reconstruction process $u_k\mapsto i_{k}^{(1)}u_k$.

Our discretization of the sea ice model is \emph{non-consisting}, which means that the discrete problem is formulated via a modified variational formulation $A_s(u_k,\phi_k)=0$ with $A_s\neq A$, see Section~\ref{sec:part} and in particular~(\ref{varsplit}) where we define the form $B_s(U)(\Phi)$. One consequence of non-conformity is the violation of Galerkin orthogonality. For the partitioned solution approach, only the following disturbed relation holds
\[
A(u,\phi_k)-A_s(u_k,\phi_k)= - A(u_k,\phi_k). 
\]
This non-consistency gives rise to the splitting terms of the error
estimator denoted by $\beta$ in~(\ref{Lsplit}).


\end{document}